\documentclass{kluwer}
\usepackage{graphicx}
\usepackage{mathrsfs}

\begin{document}
\begin{article}
\begin{opening}
\title{{\bf A Large Closed Queueing Network Containing Two Types of
Node and Multiple Customer Classes: One Bottleneck Station}}

\author{\surname{{\sc Vyacheslav M. Abramov}}\email{vyachesl@zahav.net.il}}
\institute{The Sackler Faculty of Exact Sciences, School of
Mathematical Sciences, Tel Aviv University, 69978 Tel Aviv,
Israel}




\runningtitle{Closed queueing network with two types of node}
\runningauthor{VYACHESLAV ABRAMOV}



\begin{abstract}
\noindent The paper studies a closed queueing network containing
two types of node. The first type (server station) is an infinite
server queueing system, and the second type (client station) is a
single server queueing system with autonomous service, i.e. every
client station serves customers (units) only at random instants
generated by strictly stationary and ergodic sequence of random
variables. It is assumed that there are $r$ server stations. At
the initial time moment all units are distributed in the server
stations, and the $i$th server station contains $N_i$ units,
$i=1,2,...,r$, where all the values $N_i$ are large numbers of the
same order. The total number of client stations is equal to $k$.
The expected times between departures in the client stations are
small values of the order $O(N^{-1})$~ $(N=N_1+N_2+...+N_r)$.
After service completion in the $i$th server station a unit is
transmitted to the $j$th client station with probability
$p_{i,j}$~($j=1,2,...,k$), and being served in the $j$th client
station the unit returns to the $i$th server station. Under the
assumption that only one of the client stations is a bottleneck
node, i.e. the expected number of arrivals per time unit to the
node is greater than the expected number of departures from that
node, the paper derives the representation for non-stationary
queue-length distributions in non-bottleneck client stations.
\end{abstract}
\end{opening}


\section{Introduction}
\indent

\subsection{Description of the model and motivation}

We consider a closed queueing network containing two types of
node. The first type (server station) is an infinite server
queueing system with identical servers. There are $r$ server
stations. The second type (client station) is a single server
queueing system with autonomous service, where customers (units)
are served only at random instants generated by strictly
stationary and ergodic sequence of random variables.
\smallskip

Queueing systems with autonomous service were introduced and
originally studied by Borovkov [4], [5]. Below we recall the
definition of the version of queueing system with autonomous
service in the case when both arrival of customers and their
service are ordinary (i.e. not batch), and all processes (arrival,
departure etc.) start at 0.

Let $E(t)$ and $S(t)$ be right continuous, having the left-limits,
point processes, defined for all $t\ge 0$. Let $E(t)$ describe an
{\it arrival} process, and let $S(t)$ describe a {\it departure}
process. We say that the service mechanism is {\it autonomous} if
the {\it queue-length process} $Q(t)$ is defined by the equation
$$
Q(t)=E(t)-\int_0^t{\bf I}\{Q(s-)>0\}\mbox{d}S(s),
$$
where ${\bf I}(A)$ denotes an indicator of event $A$, and the
integral is understood in the sense of the Lebesgue-Stieltjes
integral. For a more general definition of queueing systems with
autonomous service, when arrivals and departures occur by batches
see Borovkov [4], [5].

\smallskip

Along with $r$ server stations there are $k$ client stations. The
departure instants in the $j$th client station ($j=1,2,...,k$) are
denoted $\xi_{j,1}$, $\xi_{j,1}+\xi_{j,2}$,
$\xi_{j,1}+\xi_{j,2}+\xi_{j,3}$,\ldots where, as it was mentioned
above,

\medskip

    $\bullet$ each sequence $\{\xi_{j,1}, \xi_{j,2},...\}$ forms a strictly
stationary and ergodic sequence of random variables.

\medskip

The closed network contains $N$ units. At the initial time moment
they are distributed in the server stations, and the $i$th server
station contains $N_i$ units ($N_1+N_2+...+N_r=N$). The service
time of each unit of the $i$th server station is exponentially
distributed random variable with the expectation $\lambda_i$.
After a service completion at the $i$th server station a unit is
transmitted to the client station $j$ with probability $p_{i,j}\ge
0$, $\sum_{j=1}^kp_{i,j}=1$, and then, after the service
completion at the $j$th client station, the unit returns to the
$i$th server station. Thus, the model of the network considered
here is a model with {\it multiple customer classes} defined by
$r$ server and $k$ client stations.

Denote $\lambda_{i,j}=\lambda_ip_{i,j}$ and
$(\mu_jN)^{-1}=\mathbb{E}\xi_{j,1}$. Then, the input rate to the
$j$th client station is
$$
\sum_{i=1}^r\lambda_{i,j}N_i,
$$
and the traffic intensity
$$
\varrho_j(N)=\frac{1}{\mu_jN}\sum_{i=1}^r\lambda_{i,j}N_i.
$$

  In the following it will be assumed that

\medskip

  $\bullet$ the series
parameter $N$ increases to infinity, and, as $N\to\infty$, each
fraction $N_i/N$ converges to some positive number $\alpha_i$.

\medskip

Then client station $j$ is called {\it non-bottleneck station} if,
as $N\to\infty$, the limiting value
$\varrho_j=\lim_{N\to\infty}\varrho_j(N)$ is less than 1.
Otherwise, the $j$th client station is called {\it bottleneck}
station.

     It is assumed in the paper that

\medskip

    $\bullet$     the first $k-1$ client stations are non-bottleneck
stations, while the $k$th client station is a bottleneck station.

\medskip

It was assumed above that

\medskip
    $\bullet$ the probabilities $p_{i,j}$ satisfy the following
    two conditions: $p_{i,j}\ge 0$ and $\sum_{j=1}^kp_{i,j}=1$.

\medskip

These two conditions define a general class of connections between
the client and server stations. For example, if $p_{i,j}=0$ for
some indexes $i$, $j$, then there is no connection between the
$i$th server station and the $j$th client station. Keeping in mind
that according to the convention that the $k$th client station is
a bottleneck station, the behavior of queues in the non-bottleneck
client stations essentially depends on topology of the network,
that is on existing connections between the server and client
stations. In order to clarify this, consider two different
topologies of the network containing two server and four client
stations. These two topologies are shown in Figures 1.1 and 1.2
below.

\smallskip

In the network described in Figure 1.1 the probabilities
$p_{1,3}$, $p_{1,4}$, $p_{2,1}$ and $p_{2,2}$ are equal to 0,
while the probabilities $p_{1,1}$, $p_{1,2}$, $p_{2,3}$ and
$p_{2,4}$ are positive. Then four client stations are separated
into two nonintersecting groups, and the network is decomposed
into two subnetworks forming the simplest tree. The topology of
the network looks as follows.

 \scalebox{0.75}{\includegraphics{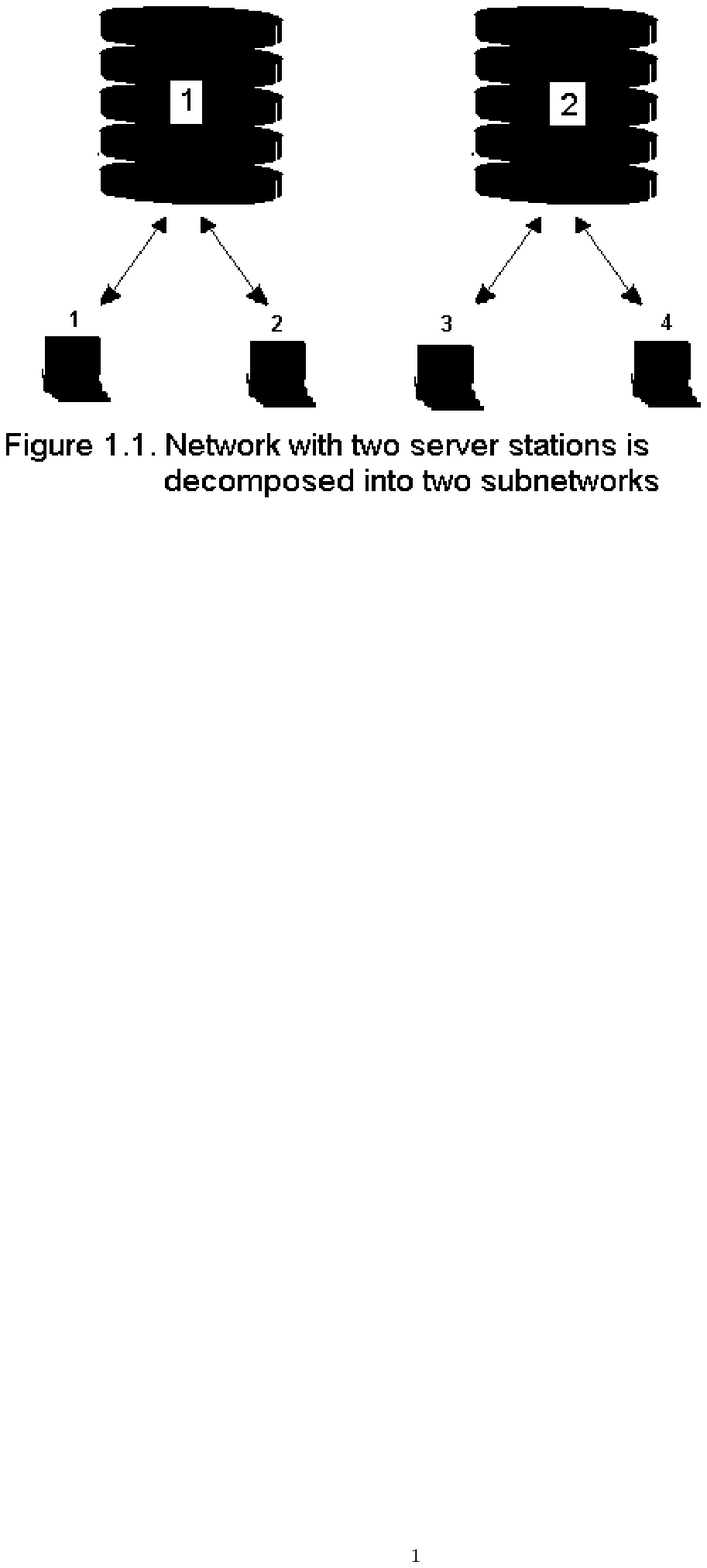}}

Another example of the network is described in Figure 1.2. In this
network all the probabilities $p_{i,j}$ are strictly positive, and
then the network can not be decomposed into subnetworks, and it
forms a net. The topology of the network looks as follows.

 \scalebox{0.75}{\includegraphics{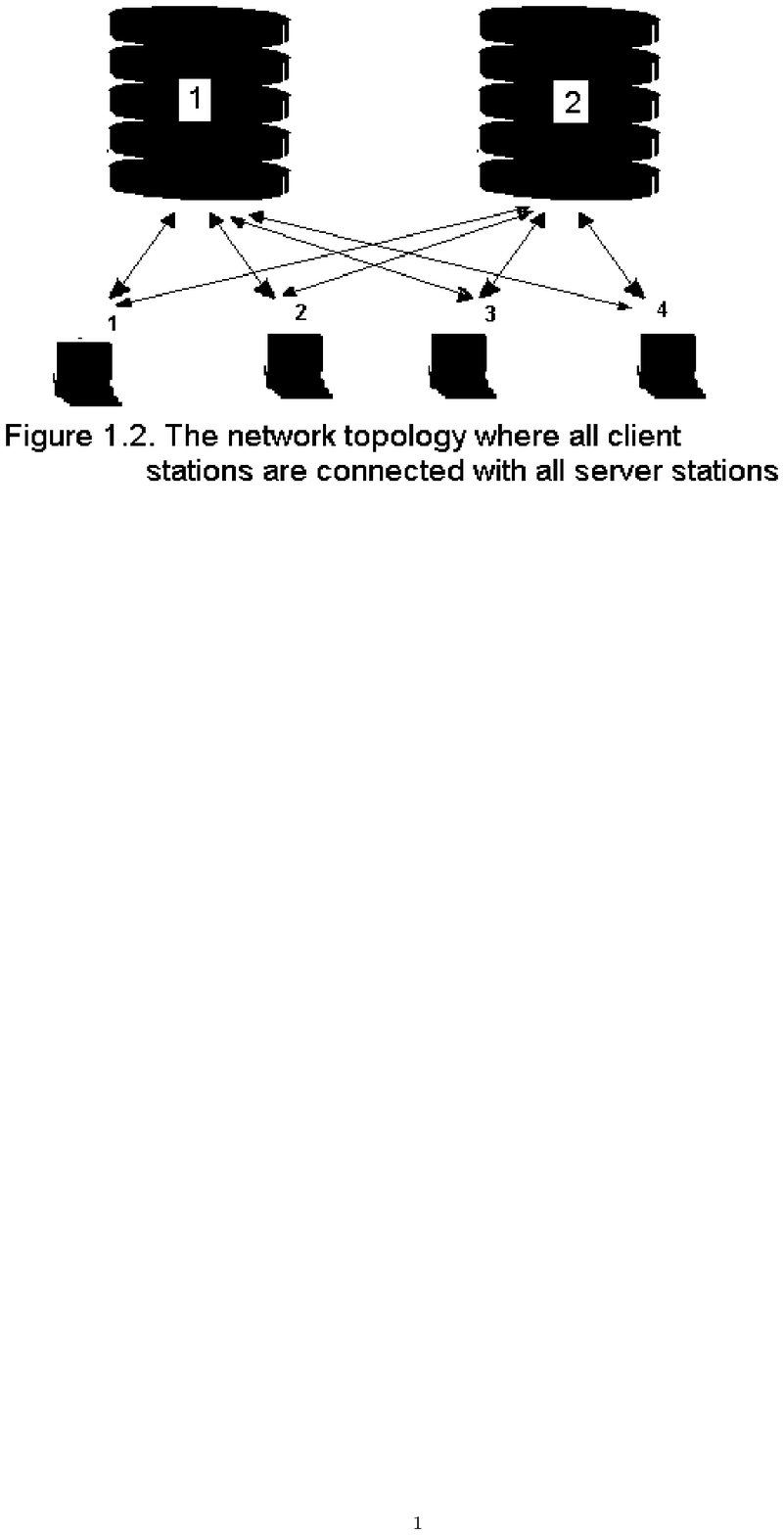}}

It is clear that the case of the network described in Figure 1.1
is artificial rather than realistic: the considered network is
explicitly separated into two subnetworks. The analysis of that
network is reduced to the same analysis of the network with one
server. One of subnetwork has a bottleneck client station, other
has non-bottleneck client stations only. The analysis of the
network with only one server station (hub) and a number of client
(satellite) stations has been done in Abramov [1], and earlier, in
the case of Markovian network in Kogan and Liptser [16].
\smallskip

The network topology described in Figure 1.2 is a topology where
there are connections between all client stations and all server
stations. The behavior of this network, as $N\to\infty$, also has
simple intuitive explanation and can be reduced to the case of the
network studied in the earlier paper of Abramov [1]. In the
following we will return to the explanation of that behavior.

\smallskip

Now, speaking about the general network topology it is worth
noting the following. The specific feature of the network, having
a hub and several satellite stations, is that if there is a
bottleneck node, then the behavior of queues in {\it all}
non-bottleneck stations depends on the behavior of the queue in
the bottleneck station. In the limiting case as $N\to\infty$ the
dependence vanishes. However, the limiting non-stationary
queue-length distributions in {\it all} non-bottleneck satellite
stations are dependent of time $t$ (see Abramov [1]). The behavior
of the network containing several server and client stations is
more complicated. For details see Sections 1.3 and 1.4.

\smallskip

There is a large number of different applications for the queueing
network with two types of node described above. Below we provide
one of the possible motivation as a database distributed in $r$
server stations. Then by unit we mean records of data or data
units that can be called by users from the client stations. The
action of user is to call and update a data unit. Being called by
one client station the data unit becomes blocked and not
accessible from other ones. After completion of the processing of
the data unit the system unblocks it, and again the data unit
becomes accessible for new actions. Note, that our assumption that
the service mechanism is autonomous is one of possible
generalizations of Markovian networks. Autonomous service
mechanism is of interest for technologies of computer systems,
where from time to time the system automatically looks up the
queue in order to provide then the service for its units
(messages, queries) waiting for their processing.

\subsection{The history of question and the goal of the paper}

The history of subject is very rich. There is a large number of
papers in the queueing literature explicitly or implicitly related
to the subject. A large number of papers study state-dependent and
time-dependent queueing models, providing different approximations
including diffusion and fluid approximations (e.g. Krylov and
Liptser [22], Mandelbaum and Massey [26], Mandelbaum {\it et al}
[27], Mandelbaum and Pats [28], [29], Chen and Mandelbaum [7],
[8], Chao [6], Williams [35], [36], [37] and others). A
state-dependent queueing model, where the arrival and service
rates depend on the current workload of the system, can be
considered as a model for separated client station of computer
network. The detailed asymptotic analysis of such systems by
analytic methods is given in a number of papers of Knessl {\it et
al} [10] - [14] etc. A number of papers study one server and one
client stations only. Containing two types of node these system
are closer to the considered system. Such systems have been
studied in Krichagina {\it et al} [20], Krichagina and Puhalskii
[21], Liptser [23] and other papers. Different models, closed to
the considered model, have been also studied in Whitt [34], Knessl
and Tier [15], Kogan {\it et al} [17], [18], Reiman and Simon [31]
and others.

\smallskip

Nevertheless, speaking about the history of subject we restrict
our attention only with a small number of papers, having a close
relation to the considered paper in some aspects as description of
model, main goal of research and mathematical methods. Only the
papers by Kogan and Liptser [16], Abramov [1], [2] will be
discussed here. These three papers form a chronological chain of
papers providing a comprehensive analysis of the queueing networks
with one server station (hub) and a number of client stations
(satellite nodes), where one of the satellite nodes is a
bottleneck station. All these networks depend on the large
parameter $N$ - the number of tasks in the server, and, as
$N\to\infty$, the papers study the limiting non-stationary
queue-length distribution in non-bottleneck satellite stations.
The papers also provide the diffusion and fluid approximations
(given under appropriate assumptions) for the queue-length in the
bottleneck satellite station. The method of analysis in these
papers is the theory of martingales, where in two of these three
papers, Kogan and Liptser [16] and Abramov [1], the dominated
method is the stochastic calculus, while the main method in
Abramov [2] is the up- and down-crossings approach and the
martingale techniques in discrete time. The part of the paper of
Abramov [2], related to the bottleneck analysis, uses the
stochastic calculus as well, but its application is related to the
special examples rather than to the theory. Below we briefly
discuss the part of the results related to the limiting
non-stationary distributions in the non-bottleneck nodes.

\smallskip

Kogan and Liptser [16] study a closed Markovian queueing network
with an $M/M/\infty$ hub, a large number $N$ of customers (tasks)
which at the initial time moment all distributed in the server,
and $k$ different $M/M/1$ satellite nodes, assuming that the $k$th
satellite node operates as a bottleneck node. It was shown that,
as $N\to\infty$, the limiting non-stationary queue-length
distribution in the non-bottleneck node is a geometric
distribution with parameter depending on time.

\smallskip

Abramov [1] develops the model considered in the paper of Kogan
and Liptser [16], assuming that the service mechanism in satellite
nodes is autonomous, and the sequence of intervals between service
completions there forms a strictly stationary and ergodic sequence
of random variables. It is established the expression, where the
limiting non-stationary queue-length distribution in
non-bottleneck node in time $t$ is expressed via limiting
non-stationary queue-length distribution immediately before the
last departure of a unit before time $t$. The obtained
representation enables us to conclude that for any time $t>0$ the
queue-length in the bottleneck node has the same order as that in
the Markovian variant of network, having the same traffic
intensities in the satellite nodes.

\smallskip

Abramov [2] develops the Markovian model of Kogan and Liptser [16]
as follows. Whereas in the model of Kogan and Liptser [16] the hub
is an infinite-server queueing system, Abramov [2] studies the
model where service times in the hub are generally distributed
random variables, depending on the number of customers (tasks)
residing there, but as in the model of Kogan and Liptser [16], the
service times in the satellite stations are assumed to be
exponentially distributed as well. More precisely, it is assumed
the following: if immediately before a service of a sequential
customer the queue-length at the hub is equal to $K\le N$, then
the probability distribution function is $G_K(Kx)$,
$g_K^{-1}=\int_0^\infty x\mbox{d}G_K(x)<\infty$ and, as
$N\to\infty$, the sequence of probability distributions
$\{G_N(x)\}$ converges weakly to $G(x)$ with $g^{-1}=\int_0^\infty
x\mbox{d}G(x)<\infty$, and $G(0+)=0$. Along with these assumptions
it is required the stochastic order relations between two neighbor
distribution functions $G_K(Kx)$ and $G_{K+1}(Kx+x)$. It is
assumed that $G_K(Kx)\le G_{K+1}(Kx+x)$ for all $x\ge 0$. The
sense of this order relation is intuitively clear: {\it a rate of
service time at the hub increases, as a queue-length increases
there.} Note, that this assumption is automatically implied in the
special case of
$G_1(x)=G_2(x)=...=G_N(x)=G(x)=1-\mbox{e}^{-\lambda x}$, leading
to the network considered by Kogan and Liptser [16].

\smallskip

Studying a more general model than in the paper of Abramov [1], we
are not going to provide a comprehensive analysis of the network
as it is done in the three abovementioned papers of Kogan and
Liptser [16] and Abramov [1], [2]. Namely, we are not going to
study diffusion and fluid approximations for a bottleneck node in
the corresponding cases of {\it moderate} and {\it heavy} usage
regimes respectively (for the definition see e.g. Kogan and
Liptser [16]). The reason for this is the following. During the
last decades the diffusion and fluid approximations have been
intensively studied in a large number of works related to quite
general class of stochastic models. Although, to our best
knowledge, the model considered here is not covered, the behavior
of the bottleneck client station under appropriate heavy traffic
conditions remains the same as of the earlier models considered in
the abovementioned papers. In other words, the behavior of the
bottleneck client station is expected to be described by the
similar stochastic It\^o equations as in the paper of Kogan and
Liptser [16] or Abramov [1]. In the analysis of the multi-server
stochastic network considered here, the new effects can be
expected namely for the limiting non-stationary queue-length
distributions in the non-bottleneck stations, and therefore, it is
the main object for study of the present paper.

\subsection{The main result and the special cases}

Introduce the following notation, which will be used throughout
the paper. The queue-length process in the client station $j$
($j=1,2,...,k$) will be denoted $Q_j(t)$. Each queue-length
process $Q_j(t)$ is a function of the set of parameters $N_1$,
$N_2$,..., $N_r$, i.e. $Q_j(t)=Q_j(t, N_1, N_2,..., N_r)$. To
avoid the complicated notation, these parameters will be always
omitted. Thus, writing $\lim_{N\to\infty}\mathbb{P}\{Q_j(t)=l\}$
we mean $\lim_{N\to\infty}\mathbb{P}\{Q_j(t$, $N_1$, $N_2$, ...,
$N_r)=l\}$ or $\lim_{N\to\infty}\mathbb{P}\{Q_j(t, N)=l\}$. (All
the parameters $N_j$, $j=1,2,...,r$ depend on $N$, and, as
$N\to\infty$, all them increase to infinity as well.) We hope that
it does not cause a misunderstanding. The same note holds for the
other processes appearing in the paper.

The point processes $S_j(t)$, $j=1,2,...,k$, associated with the
random sequences $\{\xi_{j,1}$, $\xi_{j,2}$, $...\}$, are defined
as follows:
$$
S_j(t)=\sum_{l=1}^\infty{\bf I}\{\sigma_{j,l}\le t\},
$$
where ${\bf I}(A)$ denotes the indicator of event $A$ and
$$
\sigma_{j,l}=\sum_{i=1}^l\xi_{j,i},~~~l\ge 1.
$$
For any $t>0$ introduce the process
$$
S^*(t)=\inf\{s>0: S_j(s)=S_j(t)\},
$$
having a sense as the moment of the last jump of the point process
$S_j(t)$ before time $t$.

\medskip

\noindent {\bf Theorem 1.1.} Under the assumptions of the paper
for $j=1,2,..., k-1$ and for any $t>0$ we have:
$$
\lim_{N\to\infty}\mathbb{P}\{Q_j[S_j^*(t)]=0\}=1-\rho_j(t),
$$
$$
\lim_{N\to\infty}\int_0^t\rho_j(s)\mathbb{
P}\{Q_j(s)=l\}\mbox{d}s
$$
$$
=\lim_{N\to\infty}\int_0^t\mathbb{
P}\{Q_j[S_j^*(s)]=l+1\}\mbox{d}s,
$$
$$
l=0,1,...,
$$
where
$$
\rho_j(t)=\varrho_j\Big[1-q(t)\sum_{i\in\mathscr{
I}_j\bigcap\mathscr{I}_k}\beta_{i,j}\beta_{i,k}\Big], \eqno (1.1)
$$
$$
q(t)=\Big(1-\frac{1}{\varrho_k}\Big)(1-\mbox{e}^{-\varrho_k\mu_kt}),
$$
$$
\beta_{i,j}=\frac{\lambda_{i,j}\alpha_i}{\varrho_j},~~~i=1,2,...,r,
$$
and $\mathscr{I}_j$ is the set of indexes $i=1,2,...,r$ where
$\lambda_{i,j}>0$.

\medskip

Some special cases associated with Theorem 1.1 are given below.
For Markovian network we have the following:

\medskip
\noindent {\bf Corollary 1.2.} If the point processes $S_j(t)$,
$j=1,2,...,k$, all are mutually independent Poisson processes,
then for all $j=1,2,..., k-1$
$$
\lim_{N\to\infty}\mathbb{
P}\{Q_j(t)=l\}=\big[1-\rho_j(t)\big]\big[\rho_j(t)\big]^l,
$$
$$
l=0,1,2,...,
$$
where $\rho_j(t)$ are given by (1.1).

\medskip

In the case when there is only one server station, we obtain
Theorem 3 of Abramov [1].

\medskip
\noindent {\bf Corollary 1.3.} [Abramov [1].] If $i=1$ then for
all $j=1,2,..., k-1$ and for any $t>0$
$$
\lim_{N\to\infty}\mathbb{P}\{Q_j[S_j^*(t)]=0\}=1-\varrho_j[1-q(t)],
$$
$$
\lim_{N\to\infty}\int_0^t\varrho_j[1-q(s)]\mathbb{
P}\{Q_j(s)=l\}\mbox{d}s
$$
$$
=\lim_{N\to\infty}\int_0^t\mathbb{
P}\{Q_j[S_j^*(s)]=l+1\}\mbox{d}s,
$$
$$
l=0,1,...~.
$$

The following two corollaries considers the cases when the
bottleneck and non-bottleneck client stations respectively have or
does not have common server stations. In the first case we have

\medskip

\noindent {\bf Corollary 1.4.} Assume that the client stations $j$
and $k$ are connected with the common server stations, i.e. both
$p_{i,k}$ and $p_{i,j}$ are positive for the same set of indexes
$i$. Then for all $j=1,2,..., k-1$ and for any $t>0$
$$
\lim_{N\to\infty}\mathbb{P}\{Q_j[S_j^*(t)]=0\}=1-\varrho_j[1-q(t)],
$$
$$
\lim_{N\to\infty}\int_0^t\varrho_j[1-q(s)]\mathbb{
P}\{Q_j(s)=l\}\mbox{d}s
$$
$$
=\lim_{N\to\infty}\int_0^t\mathbb{
P}\{Q_j[S_j^*(s)]=l+1\}\mbox{d}s,
$$
$$
l=0,1,...~.
$$

In the second case we have
\medskip

\noindent{\bf Corollary 1.5.} If for the client stations $j$ and
$k$ there is no common server station, i.e. the subset
$\mathscr{I}_j\bigcap\mathscr{I}_k=\emptyset$, then for all
$j=1,2,..., k-1$ and for any $t>0$
$$
\lim_{N\to\infty}\mathbb{P}\{Q_j[S_j^*(t)]=0\}=1-\varrho_j, \eqno
(1.2)
$$
$$
\lim_{N\to\infty}\int_0^t\mathbb{
P}\{Q_j(s)=l\}\mbox{d}s
$$
$$
=\frac{1}{\varrho_j}\lim_{N\to\infty}\int_0^t\mathbb{
P}\{Q_j[S_j^*(s)]=l+1\}\mbox{d}s,\eqno (1.3)
$$
$$
l=0,1,...~,
$$
i.e. limiting non-stationary queue-length distribution in that
client station is independent of time, and therefore it coincides
with the limiting stationary queue-length distribution.

\medskip

Finally, we have

\medskip
\noindent{\bf Corollary 1.6.} If $\varrho_k=1$, then the limiting
non-stationary queue-length distribution in all client stations
$j$, $1\le j\le k-1$, is independent of time and given by (1.2)
and (1.3).

\subsection{Discussion of the main result on simple examples}

To discuss the main result of the paper and the corollaries we
give a number of simple examples of network topologies of the
network containing two server stations and four client stations.
The two simple examples for topologies of that network have been
considered above in Figures 1.1 and 1.2. The case of the network
topology in Figure 1.1 is intuitively clear and can be described
without any analysis by applying the known results on the network
with one server station.

Let us consider the network topology in Figure 1.2. There are four
client and two server stations, and the fourth client station is a
bottleneck node. The two server stations of that network are
common for all client stations. Therefore, according to Corollary
1.4, the limiting non-stationary queue-length distribution of the
non-bottleneck client station is as in the network with a single
server station, illustrated in Figure 1.3.

 \scalebox{0.75}{\includegraphics{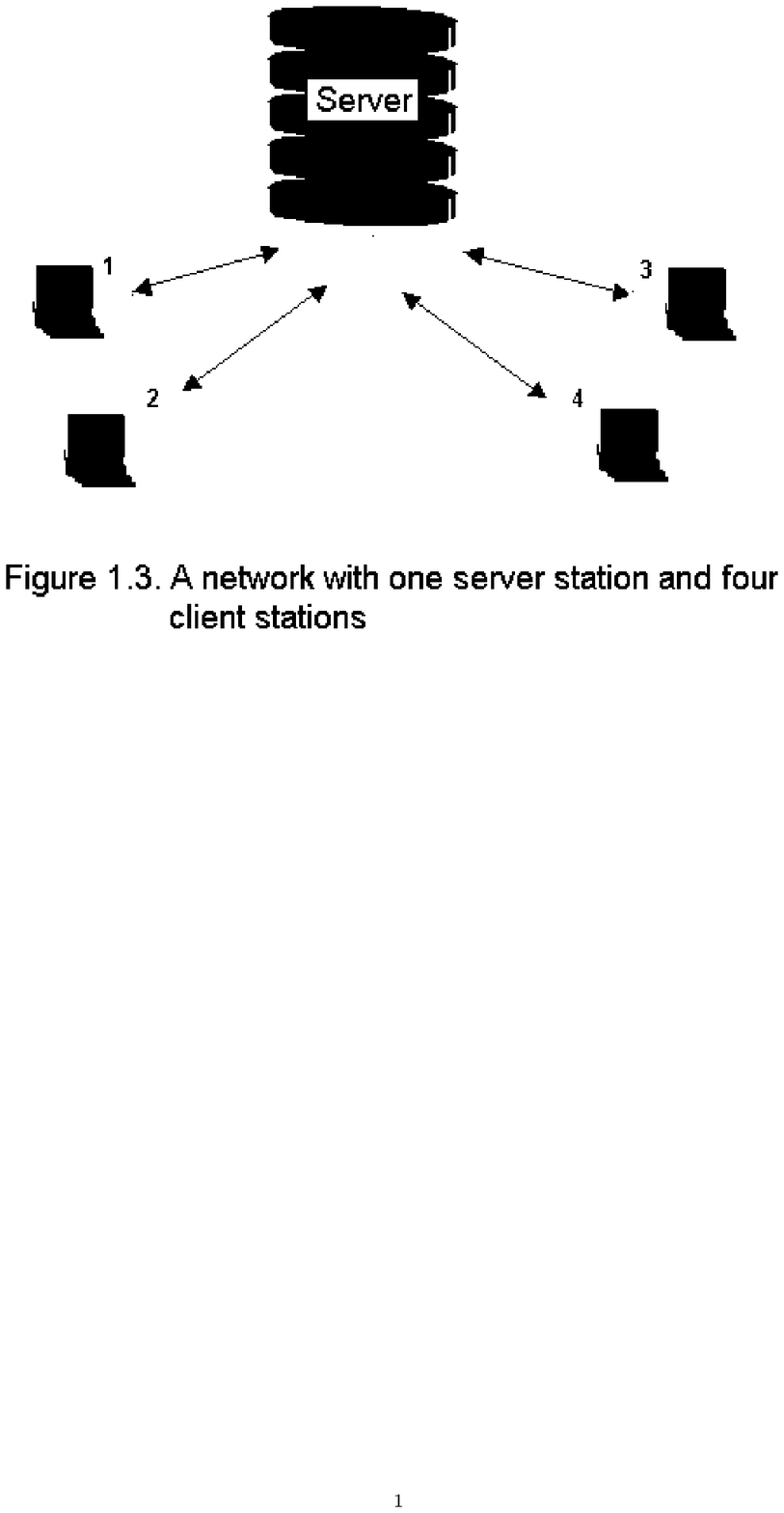}}

Thus, we join two server stations to one common server station,
assuming that the traffic parameters in the both networks are
equivalent. The intuitive explanation of this case is the
following. As $N$ large, all outputs from server stations are
close to Poisson processes. The joining of the processes outgoing
from the server stations is close to Poisson process as well.
Next, the input processes to the client stations are the thinning
of the processes outgoing from the server stations, and they are
also close to the corresponding Poisson processes. Again, the
joining of these thinning processes corresponding to the server
stations leads to the processes closed to Poissonian.

Let us now consider a new network topologies of two server and
four client stations, as it is shown in Figure 1.4. Both the first
and second server stations have connection with three client
stations. The server station 1 is not connected with a bottleneck
node. One of connections of the server station 2 is the bottleneck
client station. (According to convention we always assume that the
bottleneck node is the client station 4.)

 \scalebox{0.75}{\includegraphics{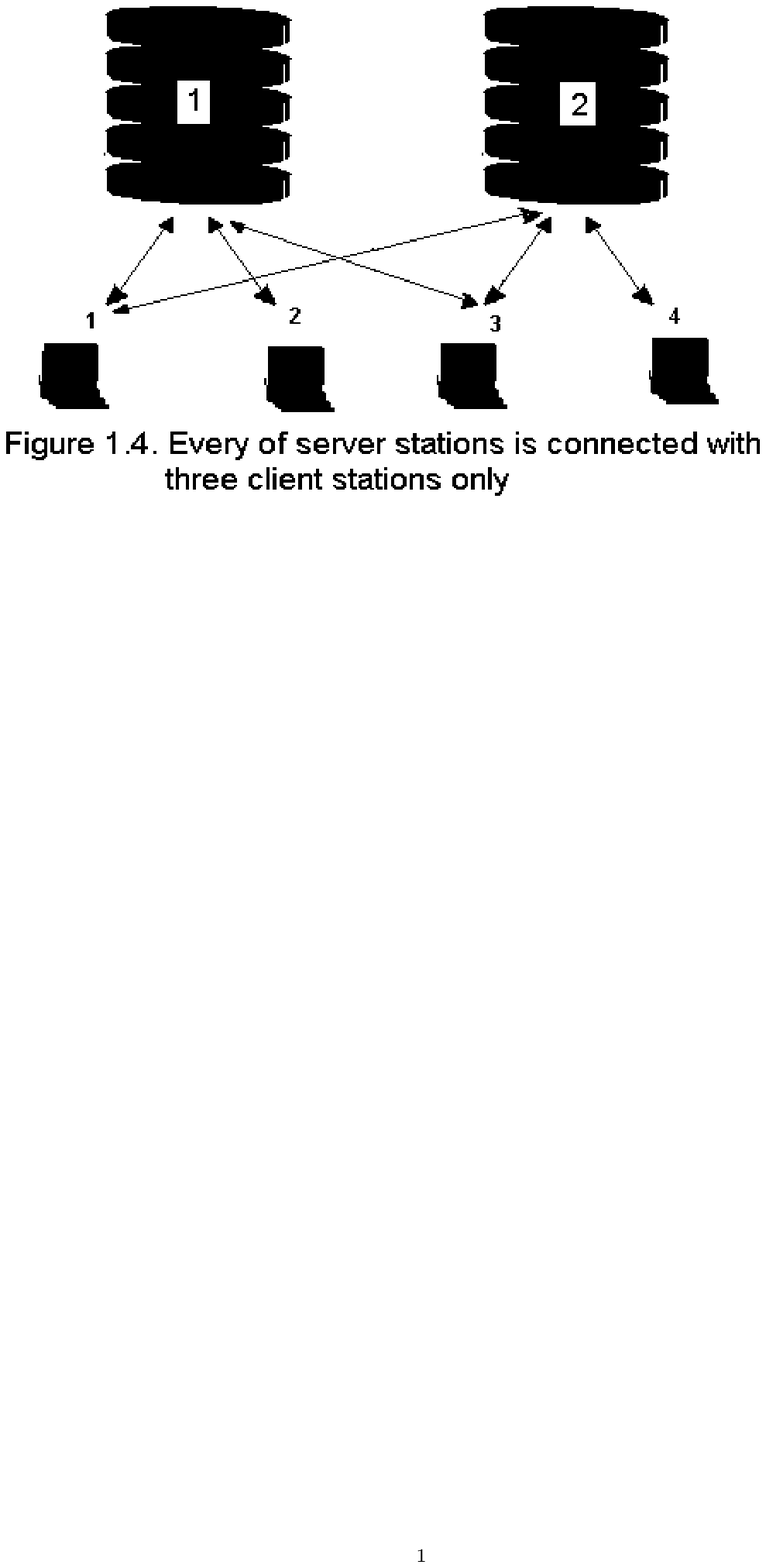}}

Then, according to Corollary 1.5, the limiting non-stationary
queue-length distribution in the client station 2 is independent
of time for any $t>0$ and coincides with the limiting stationary
queue-length distribution as $t\to\infty$. The client stations 1
and 3 depend on bottleneck station 4. There is only one server
station connected with this bottleneck station, and the client
stations 1, 3 and 4 are connected with the common server station
2. Therefore, the subnetwork consisting of the server station 2
and the client stations 1, 3 and 4 can be considered as a network
with a single server station and three client stations, and the
limiting non-stationary queue-length distribution can be
calculated by Corollary 1.4. The intuitive explanation for the
case of the abovementioned subnetwork of the server station 2 and
the client station 1, 3 and 4 is not the same as in the case of
the network in Figure 1.2, since input rate to the client stations
1 and 3 includes both streams from the server stations 1 and 2
while the input stream to the client station 4 arrives from the
server station 2 only.

Our last example is the case, where the network can not be reduced
to the more simple cases as it was above. Let us consider the
network in Figure 1.5.

 \scalebox{0.75}{\includegraphics{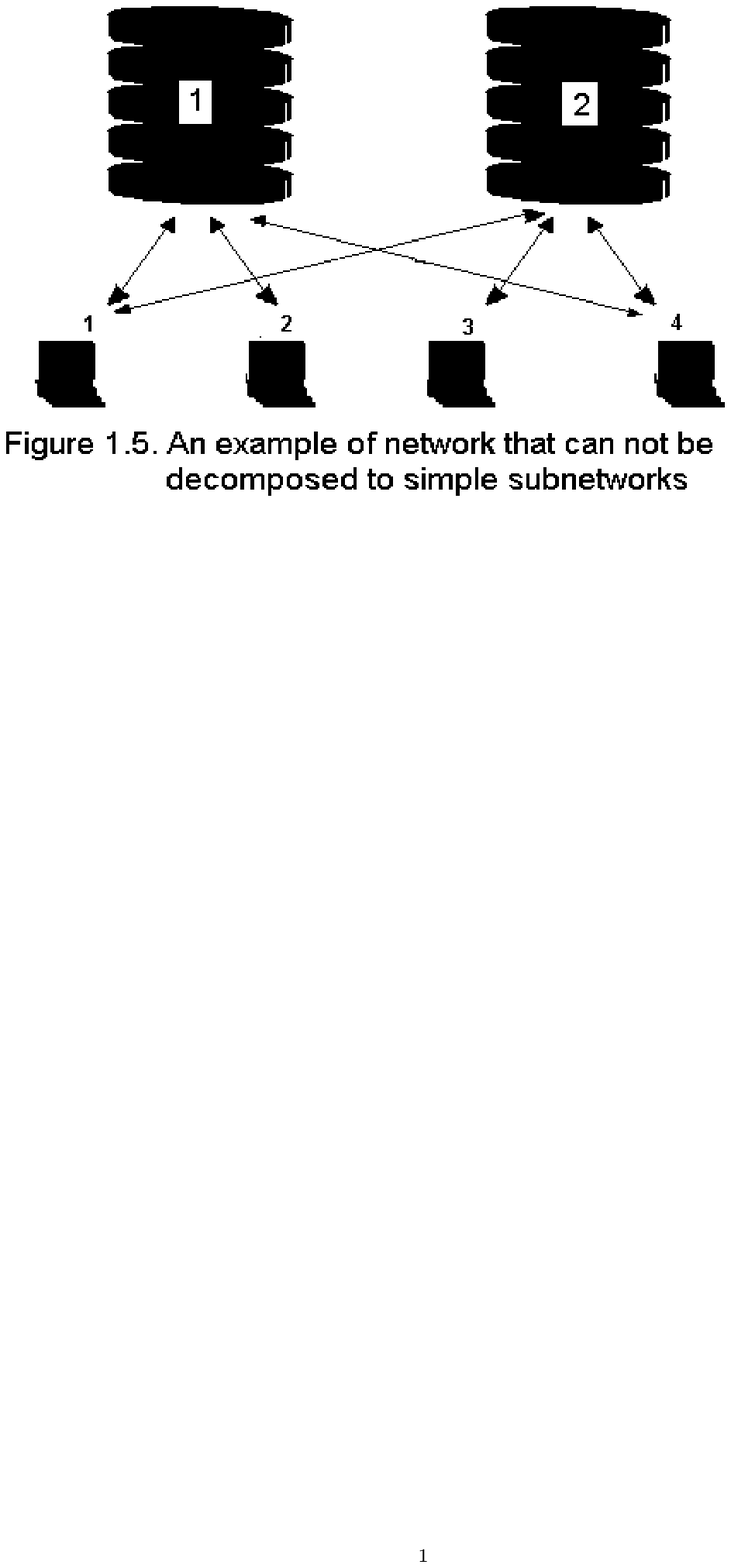}}

The bottleneck node - client station 4, receives units from the
both server stations. Every of client stations is connected only
with one server station. Therefore, neither Corollary 1.4 nor
Corollary 1.5 can be applied. Then, as $N\to\infty$, the limiting
non-stationary queue-length distributions for the client stations
1, 2 and 3 are determined from the main result, Theorem 1.1.

\subsection{The structure of the paper}

The paper is structured into 5 sections. Section 1 is an
introduction, where description of the model, review of the
literature, motivation, main result and its discussion are
provided. Section 2 derives the equations for the queue-length
processes in the client sections and reduces the problem to the
Skorokhod reflection principle. Section 3 studies asymptotic
properties of the normalized queue-length processes in the client
stations. It deduces a system of equations of the normalized
queue-lengths, permitting us to describe a dynamics of the
normalized queue-lengths in the client stations. Section 4 proves
the stability of the queue-length process in non-bottleneck client
stations, i.e. existence of functionals associated with the
limiting non-stationary queue-length distributions in those client
stations. The main result of this work is proved in Section 5.

\section{Queue-length processes in the client stations and
reduction to the Skorokhod problem}

\noindent Consider the client station $j$,~ $j=1,2,...,k$. Recall
that $Q_{j}(t)$ denote a queue-length at time $t$. According to
the convention, $Q_j(0)=0$, and for positive time instant $t>0$ we
have:
$$
Q_j(t)=A_j(t)-D_j(t), \eqno (2.1)
$$
where $A_j(t)$ is the {\it arrival} process to client station $j$
and $D_j(t)$ is the {\it departure} process from client station
$j$.

The departure process is described by the equation
$$
D_j(t)=\int_0^t{\bf I}\{Q_j(s-)>0\}\mbox{d}S_j(s)
$$
$$
~~~~~~~~~~~~~~~~~~=S_{j}(t)-\int_0^t{\bf
I}\{Q_j(s-)=0\}\mbox{d}S_j(s),\eqno (2.2)
$$
where the point process $S_j(t)$ was defined above in Section 1.4.

The description of the arrival process $A_j(t)$ is the following.
Let $\mathscr{I}_j$ be the set of indexes $i$ where
$\lambda_{i,j}>0$, $i=1,2,...,r$ (see also Section 1.4). Then we
have the representation
$$
A_j(t)=\sum_{i\in\mathscr{I}_j}A_{i,j}(t), \eqno (2.3)
$$
where $A_{i,j}(t)$ denotes the arrival process from the server
station $i$ to the client station $j$.

In order to derive the explicit representations for the processes
$A_{i,j}(t)$ let us introduce the following notation. Let
$\mathscr{J}_i$ be the set of indexes $j$ where $\lambda_{i,j}>0$,
and let $Q_{i,j}(t)$ denote the number of units in the queue at
the $j$th client station arriving from the $i$th server station.
Introduce also $\{\pi_{i,j,v}(t)\}$, $v=1,2,...,N_i$, a collection
of independent Poisson processes with rate $\lambda_{i,j}$. Then
we have
$$
A_{i,j}(t)=\int_0^t\sum_{v=1}^{N_i}{\bf
I}\Big\{N_i-\sum_{l\in\mathscr{J}_i}Q_{i,l}(s-)\ge
u\Big\}\mbox{d}\pi_{i,j,v}(s), \eqno (2.4)
$$
and
$$
Q_j(t)=\sum_{i\in\mathscr{I}_j}Q_{i,j}(t). \eqno (2.5)
$$

From (2.1), (2.2) we have the following representation for the
queue-length process $Q_j(t)$
$$
Q_j(t)=A_j(t)-S_j(t)+\int_0^t{\bf I}\{Q_j(s-)=0\}\mbox{d}S_j(t),
\eqno (2.6)
$$
which implies that $Q_j(t)$ is the normal reflection of the
process
$$X_j(t)=A_j(t)-S_j(t),~~~~~ X_j(0)=0
\eqno (2.7)
$$
at zero. More accurately, $Q_j(t)$ is a non-negative solution of
the Skorokhod problem (e.g. Skorokhod [32], Tanaka [33], Anulova
and Liptser [3] of the normal reflection of the process $X_j(t)$
at zero (for the detailed arguments see Kogan and Liptser [16]).
According to Skorokhod problem
$$
\varphi_j(t)=\int_0^t{\bf
I}\{Q_j(s-)=0\}\mbox{d}S_j(s)=-\inf_{s\le t}X_j(s). \eqno (2.8)
$$
In the following we will use the notation
$$
\Psi_t(X)=-\inf_{s\le t}X(s) \eqno (2.9)
$$
for any c\'adl\'ag function $X(t)$, $t\ge 0$, with $X(0)=0$. Then
we have
$$
Q_j(t)=X_j(t)+\Psi_t(X_j). \eqno (2.10)
$$
We will use also the notation
$$\Phi_t(X)=X(t)+\Psi_t(X) \eqno (2.11)
$$
for any c\'adl\'ag function $X(t)$, $t\ge 0$, with $X(0)=0$. Then
(2.11) is rewritten as $Q_j(t)=\Phi_t(X_j)$.

\smallskip

Take into account that the processes $A_{i,j}(t)$, $S_j(t)$,
$j=1,2,...,k$, $i\in\mathscr{I}_j$,  all are the semimartingales
adapted with respect to filtration $\mathscr{F}_t$ given on
stochastic basis $\{\Omega, \mathscr{F}, {\bf
F}=(\mathscr{F}_t)_{t\ge 0}, \mathbb{P}\}$. In the following the
compensators of these processes are pointed out by hat. For
example, $\widehat S_j(t)$ is a compensator corresponding to the
semimartingale $S_j(t)$ in the Doob-Meyer decomposition (i.e.
$S_j(t)=\widehat S_j(t)+M_{S_j}(t)$). Taking into account that the
processes $S_j(t)-\widehat S_j(t)$, $A_{i,j}(t)-\widehat
A_{i,j}(t)$, $j=1,2,...,k$, $i\in\mathscr{I}_j$, all are the local
square integrable martingales (see Liptser and Shiryayev [24],
Chapter 18) we obtain
$$
X_j(t)=\widehat A_j(t)-\widehat S_j(t)+M_j(t), \eqno (2.12)
$$
where
$$
\widehat A_j(t)=\sum_{i\in\mathscr{I}_j}\widehat A_{i,j}(t) \eqno
(2.13)
$$
(see ref. (2.3)), and
$$
M_j(t)=[A_j(t)-\widehat A_j(t)]-[S_j(t)-\widehat S_j(t)] \eqno
(2.14)
$$
is a local square integrable martingale. In addition to (2.12)
$$
\widehat A_{i,j}(t)=\int_0^t\lambda_{i,j}\Big\{N_i-
\sum_{l\in\mathscr{J}_i}Q_{i,l}(s)\Big\}\mbox{d}s \eqno (2.15)
$$
(for details see Dellacherie [9], Liptser and Shiryayev [25],
theorem 1.6.1). Introduce the random function $g_{i,j}(t)$, the
instantaneous rate of units arriving from the server station $i$
to the client station $j$ in time $t>0$
$$
g_{i,j}(t)=\lim_{\Delta\to 0}\frac{\mathbb{E}\{
A_{i,j}(t)-A_{i,j}(t-\Delta)|N_i-
\sum_{l\in\mathscr{J}_i}Q_{i,l}(t)\}}{\Delta}
$$
$$
~~~~=\lim_{\Delta\to 0}\frac{\widehat A_{i,j}(t)-\widehat
A_{i,j}(t-\Delta)}{\Delta}
$$
$$
=\lambda_{i,j}\Big\{N_i-
\sum_{l\in\mathscr{J}_i}Q_{i,l}(t)\Big\},\eqno (2.16)
$$
coinciding with the integrand of (2.15).
 Then, the sense of
$$
\frac{g_{i,j}(t)}{\sum_{l\in\mathscr{I}_j}g_{l,j}(t)} \eqno (2.17)
$$
is the fraction of the instantaneous rate of units arriving from
the server station $i$ to the client station $j$ in time $t$ with
respect to the instantaneous rate of units arriving to the client
station $j$ in time $t$.

For small $\Delta$ and $t-\Delta\ge 0$ let
$$
M_{i,j}(t)-M_{i,j}(t-\Delta)=\frac{g_{i,j}(t)}{\sum_{l\in\mathscr{
I}_j}g_{l,j}(t)}[A_{j}(t)-A_{j}(t-\Delta)-\widehat
A_{j}(t)+\widehat A_{j}(t-\Delta)]
$$
$$~~~~~~~~~~~~~~~~~~~~~~~~~~~~~~~~~~~~
-\frac{g_{i,j}(t)}{\sum_{l\in\mathscr{
I}_j}g_{l,j}(t)}[S_{j}(t)-S_{j}(t-\Delta)-\widehat
S_{j}(t)+\widehat S_{j}(t-\Delta)].\eqno (2.18)
$$
Then
$$
M_{i,j}(t)= \int_0^t\frac{g_{i,j}(s)}{\sum_{l\in\mathscr{
I}_j}g_{l,j}(s)}~\mbox{d} M_j(s) \eqno (2.19)
$$
is a local square integrable martingale. It is readily seen from
(2.19) that
$$
\sum_{i\in \mathscr{I}_j}M_{i,j}(t)=\sum_{i\in \mathscr{
I}_j}\int_0^t\frac{g_{i,j}(s)}{\sum_{l\in\mathscr{
I}_j}g_{l,j}(s)}~\mbox{d} M_j(s)
$$
$$~~~~~~~~~~~~~~~
=\int_0^t\sum_{i\in
\mathscr{I}_j}\frac{g_{i,j}(s)}{\sum_{l\in\mathscr{
I}_j}g_{l,j}(s)}~\mbox{d} M_j(s)
$$
$$
=M_j(t).~~~~~~~~~~~~~\eqno (2.20)
$$

\section{Asymptotic properties of normalized queue-length in
client stations} \noindent For $j=1,2,...,k$, let
$$
q_j(t)=\frac{1}{N}Q_j(t) \eqno (3.1)
$$
denote the normalized queue-length process, and let
$$
x_j(t)=\frac{1}{N}X_j(t) \eqno (3.2)
$$
denote the associated normalization of the process $X_j(t)$. From
(2.12) we have
$$
x_j(t)=\frac{1}{N}\widehat A_j(t)-\frac{1}{N}\widehat
S_j(t)+m_j(t), \eqno (3.3)
$$
where
$$
m_j(t)=\frac{1}{N}M_j(t) \eqno (3.4)
$$
is a local square integrable martingale. From (3.1) - (3.4) we
obtain the equation
$$
q_j(t)=\frac{1}{N}\widehat A_j(t)-\frac{1}{N}\widehat
S_j(t)+m_j(t)+\Psi_t(x_j). \eqno (3.5)
$$
Following (2.11), equation (3.5) can be rewritten in the other
form
$$
q_j(t)=\Phi_t(x_j). \eqno (3.6)
$$
Going back to relation (2.5) we can also write
$$
q_j(t)=\sum_{i\in\mathscr{I}_j}q_{i,j}(t), ~~~i=1,2,...,r, \eqno
(3.7)
$$
where
$$
q_{i,j}(t)=\frac{1}{N}Q_{i,j}(t). \eqno (3.8)
$$
In the following we will use the formalization
$$
q_{i,j}(t)=\Phi_t^i(x_j), \eqno (3.9)
$$
giving us
$$
\Phi_t(x_j)=\sum_{i\in\mathscr{I}_j}\Phi_t^i(x_j). \eqno (3.10)
$$

Next, let us introduce the normalized process
$$
x_{i,j}(t)=\frac{1}{N}\int_0^t\frac{g_{i,j}(s)}{\sum_{l\in\mathscr{
I}_j}g_{l,j}(s)}~\mbox{d}\widehat A_{j}(s)
$$
$$~~~~~~~~~~~~~~~~~~~~
-\frac{1}{N}\int_0^t\frac{g_{i,j}(s)}{\sum_{l\in\mathscr{
I}_j}g_{l,j}(s)}~\mbox{d}\widehat S_j(s)+m_{i,j}(t),\eqno (3.11)
$$
where
$$
m_{i,j}(t)=\frac{1}{N}M_{i,j}(t). \eqno (3.12)
$$
Relations (3.11), (2.15) and (2.17) yield
$$
x_{i,j}(t)=\frac{1}{N}\int_0^t\lambda_{i,j}\Big\{N_i-
\sum_{l\in\mathscr{J}_i}Q_{i,l}(s)\Big\}\mbox{d}s
$$
$$~~~~~~~~~~~~~~~~~~~~
-\frac{1}{N}\int_0^t\frac{g_{i,j}(s)}{\sum_{l\in\mathscr{
I}_j}g_{l,j}(s)}~\mbox{d}\widehat S_j(s)+m_{i,j}(t).\eqno (3.13)
$$
Taking the sum over $i$ in (3.13) we obtain
$$
x_j(t)=\frac{1}{N}\int_0^t\sum_{i\in\mathscr{I}_j}\lambda_{i,j}
\Big\{N_i-\sum_{l\in\mathscr{J}_i}Q_{i,l}(s)\Big\}\mbox{d}s
$$
$$
-\frac{1}{N}\widehat S_j(t)+m_j(t).~~~~~~~~~~~~~~~~~\eqno (3.14)
$$

Let us now study the solution of equation (3.14) as $N\to\infty$.
Notice first that
$$
\mathbb{P}-\lim_{N\to\infty}\frac{1}{N}\widehat S_j(t)=\mu_jt.
\eqno (3.15)
$$
($\mathbb{P}-\lim$ denotes the limit in probability.) This
limiting relation is proved in Abramov [1], p.30-31. Here we
briefly recall the main steps of that proof.

We have
$$
\frac{1}{N}\widehat S_j(t)=\frac{1}{N}[\widehat
S_j(t)-S_j(t)]+\frac{1}{N}S_j(t)
$$
$$
=I_1(N)+I_2(N).~~~~~\eqno (3.16)
$$
The first term of the right-hand side $I_1(N)$ is a local square
integrable martingale. As $N\to\infty$ this term vanishes in
probability. The proof of that is established by application of
the Lenglart-Rebolledo inequality (see also Konstantopoulos {\it
et al} [19]). In view of (3.16) this means
$$
\mathbb{P}-\lim_{N\to\infty}\frac{1}{N}\widehat S_j(t)=\mathbb{
P}-\lim_{N\to\infty}I_2(N). \eqno (3.17)
$$
Hence, (3.15) will be proved if we prove that
$$
\mathbb{P}-\lim_{N\to\infty}\frac{1}{N}S_j(t)= \mu_jt. \eqno
(3.18)
$$
Limiting relation (3.18) in turn follows by application of the
following lemma of Krichagina {\it et al} [20].

\medskip
\noindent {\bf Lemma 3.1.} Let $\mathscr{A}^N=(\mathscr{
A}_t^N)_{t\ge 0}$, $N\ge 1$, be a sequence of increasing right
continuous random processes with $\mathscr{A}_0^N=0$. Let
$$
\mathscr{B}_t^N=\inf\{s: \mathscr{A}_s^N>t\},~~~~~t\ge 0,
$$
where $\inf(\emptyset)=\infty$. If, for every $t$ taken from dense
set $\mathscr{S}\subset{\bf R}_+$,~ $\mathscr{B}_t^N\to at$ as
$N\to\infty$~ ($a>0$), then, as $N\to\infty$,
$$
\sup_{t\le T}\Big|\mathscr{A}_t^N-\frac{t}{a}\Big|\to 0
$$
in probability for each $T>0$.\\

Note that analogously to (3.17) and (3.18) we have
$$
\mathbb{P}-\lim_{N\to\infty}\frac{1}{N}\widehat
A_{i,j}(t)=\mathbb{ P}-\lim_{N\to\infty}\frac{1}{N}A_{i,j}(t),
\eqno (3.19)
$$
and therefore
$$
\mathbb{P}-\lim_{N\to\infty}\frac{1}{N}\widehat A_j(t)=\mathbb{
P}-\lim_{N\to\infty}\frac{1}{N}A_j(t). \eqno (3.20)
$$

Next, applying the Lenglart-Rebolledo inequality, for any positive
$\delta$ we obtain
$$
\mathbb{P}\Big\{\sup_{0\le s\le t}\big|m_j(t)\big|>\delta\Big\}
$$
$$
\le\mathbb{P}\Big\{\sup_{0\le s\le
t}\Big|\big[A_j(s)+S_j(s)\big]-\big[\widehat A_j(s)+\widehat
S_j(s)\big]\Big|>\delta N\Big\}
$$
$$
\le\frac{\varepsilon}{\delta^2}+\mathbb{P}\{\widehat
A_j(t)+\widehat S_j(t)>\varepsilon N^2\}
$$
$$~~~~~~~~
\le\frac{\varepsilon}{\delta^2}+\mathbb{P}\{\widehat
S_j(t)>\varepsilon
N^2-t\sum_{i\in\mathscr{I}_j}\lambda_{i,j}N_i\},\eqno (3.21)
$$
and because of arbitrariness of $\varepsilon>0$, $m_j(t)$ vanishes
in probability as $N\to\infty$. Analogously, it is not difficult
to conclude that, as $N\to\infty$, also $m_{i,j}(t)$ vanishes in
probability.

Next, let
$$
\mathbb{P}-\lim_{N\to\infty}x_j(t)=x_j^*(t), \eqno (3.22)
$$
and
$$
\mathbb{P}-\lim_{N\to\infty}x_{i,j}(t)=x_{i,j}^*(t). \eqno (3.23)
$$
Then we have
$$
x_{i,j}^*(t)=
\int_0^t\lambda_{i,j}\Big[\alpha_i-\sum_{l\in\mathscr{
J}_i}\Phi_s^i(x_{l}^*)\Big]\mbox{d}s
$$
$$~~~~~~~~~~~~~~~~~~~~~~
-~ \mathbb{
P}-\lim_{N\to\infty}\frac{1}{N}\int_0^t\frac{g_{i,j}(s)}{\sum_{i\in\mathscr{
I}_j}g_{i,j}(s)}~\mbox{d}\widehat S_{j}(s),\eqno (3.24)
$$
and taking the sum over $i$ in view of (3.15) we obtain
$$
x_j^*(t)=\int_0^t\sum_{i\in\mathscr{
I}_j}\lambda_{i,j}\Big[\alpha_i-\sum_{l\in\mathscr{
J}_i}\Phi_s^i(x_{l}^*)\Big]\mbox{d}s + \mu_jt. \eqno (3.25)
$$
The solution of system (3.24), (3.25) is unique since because of
the Lipschitz conditions
$$
\sup_{t\le T}\big|\Phi_t(X)-\Phi_t(Y)\big|\le 2\sup_{t\le
T}\big|X_t-Y_t\big|, \eqno (3.26)
$$
and
$$
\sup_{t\le T}\big|\Phi_t^i(X)-\Phi_t^i(Y)\big|\le 2\sup_{t\le
T}\big|\Phi_t(X)-\Phi_t(Y)\big|
$$
$$~~~~~~~~~~~~~~~~~~
\le 4 \sup_{t\le T}\big|X_t-Y_t\big|.\eqno (3.27)
$$
 For the proof of Lipschitz condition (3.26) see Kogan and
Liptser [16]. In turn, the proof of (3.27) follows easily from
(3.10) and the triangle inequality.

\smallskip

Before solving the system (3.24) and (3.25) let us first establish
some properties of these equations. From (3.25) we obtain
$$
x_j^*(t)=\int_0^t\varrho_j\mu_j\Big[1-\frac{1}{\varrho_j\mu_j}\sum_{i\in\mathscr{I}_j}\sum_{l\in\mathscr{
J}_i}\lambda_{i,j}\Phi_s^i(x_{l}^*)\Big]\mbox{d}s + \mu_jt
$$
$$
\le \int_0^t\varrho_j\mu_j\Big[1-\sum_{l\in\mathscr{
J}_i}\Phi_s(x_{l}^*)\Big]\mbox{d}s + \mu_jt.~~~~~~~~\eqno (3.28)
$$
Let us now consider the system of equations
$$
\widetilde x_j^*(t)= \int_0^t\varrho_j\mu_j\Big[1-\sum_{l=1}
^k\Phi_s(\widetilde x_{l}^*)\Big]\mbox{d}s + \mu_jt, \eqno (3.29)
$$
$$
j=1,2,...,k.
$$
The unique solution of the system (3.29) is
$$
\widetilde x_k^*(t)=q(t)=
\Big(1-\frac{1}{\varrho_k}\Big)(1-\mbox{e}^{-\varrho_k\mu_k
t}),\eqno (3.30)
$$
$$
\widetilde
x_{j}^*(t)=(\varrho_j\mu_j-\mu_j)t-\varrho_j\mu_j\int_0^tq(s)\mbox{d}s,
 \eqno (3.31)
$$
$$
j=1,2,...,k-1.
$$
Then considering the sequence of processes $\widetilde
x_{j}(t)=\widetilde x_{j}(N,t)$ satisfying the system of equation
$$
\widetilde x_j(t)=\frac{1}{N}\int_0^t\varrho_j\mu_j
\Big\{1-\sum_{l=1}^k\Phi_{t}(\widetilde x_l )\Big\}\mbox{d}s
$$
$$
-\frac{1}{N}\widehat S_j(t)+m_j(t),~~~~~~~~~~\eqno (3.32)
$$
one can conclude the following. First,
$$
x_j(t)\le\widetilde x_j(t),~~~j=1,2,...,k, \eqno (3.33)
$$
where $x_j(t)$ are given by (3.14) (see arguments in ref. (3.28)).
Second, following Abramov [1], Lemma 2, for any fixed $t>0$ and $\varepsilon>0$ we have\\
$$
\lim_{N\to\infty}\mathbb{P}\Big\{\sup_{s\le t}\big|\widetilde
x_{j}(s)-\widetilde x_j^*(s)\big|\ge\varepsilon\Big\}=0, \eqno
(3.34)
$$
$$
j=1,2,...,k.
$$
From (3.34) and Lipschitz condition (3.26) we obtain
$$
\lim_{N\to\infty}\mathbb{P}\Big\{\sup_{s\le t}\Phi_s(\widetilde
x_j)\ge\varepsilon\Big\}=0, \eqno (3.35)
$$
$$
j=1,2,...,k-1,
$$
and
$$
\lim_{N\to\infty}\mathbb{P}\Big\{\sup_{s\le t}\big|
\Phi_s(\widetilde x_k)-\widetilde
x_k(s)\big|\ge\varepsilon\Big\}=0. \eqno (3.36)
$$
It follows from (3.35) and (3.36) that
$$
\Phi_t(\widetilde x_j^*)=0,~~~j=1,2,...,k-1, \eqno (3.37)
$$
and
$$
\Phi_t(\widetilde x_k^*)=q(t), \eqno (3.38)
$$
where $q(t)$ is defined above in ref. (3.30) and in Section 1.4.
Then, from Lipschitz condition (3.27) applied to (3.25), along
with (3.37) and (3.38) we obtain
$$
\Phi_t(x_j^*)=0,~~~j=1,2,...,k-1, \eqno (3.39)
$$
and
$$
\Phi_t(x_k^*)=q(t). \eqno (3.40)
$$
Substituting (3.39) for (3.24) and (3.25) we now have the
following system ($j=1,2,...,k$; $i=1,2,...,r$)
$$
x_{i,j}^*(t)= \int_0^t\lambda_{i,j}\Big[\alpha_i-
\Phi_s^i(x_{k}^*)\Big]\mbox{d}s
$$
$$~~~~~~~~~~~~~~~~~~~~~~~~~~~~
-~~ \mathbb{
P}-\lim_{N\to\infty}\frac{1}{N}\int_0^t\frac{g_{i,j}(s)}{\sum_{l\in\mathscr{
I}_j}g_{l,j}(s)}~\mbox{d}\widehat S_{j}(s),\eqno (3.41)
$$
$$
x_j^*(t)=\int_0^t\sum_{i\in\mathscr{I}_j}\lambda_{i,j}\Big[\alpha_i-
\Phi_s^i(x_{k}^*)\Big]\mbox{d}s - \mu_jt. \eqno (3.42)
$$
Let us now find the term
$$
\mathbb{
P}-\lim_{N\to\infty}\frac{1}{N}\int_0^t\frac{g_{i,j}(s)}{\sum_{l\in\mathscr{
I}_j}g_{l,j}(s)}~\mbox{d}\widehat S_{j}(s) \eqno (3.43)
$$
of (3.41) ($j=1,2,...,k$). For this purpose let us find
$$
\mathbb{P}-\lim_{N\to\infty} \frac{g_{i,j}(t)}{\sum_{l\in\mathscr{
I}_j}g_{l,j}(t)}. \eqno (3.44)
$$

From the initial condition
$$
\mathbb{
P}-\lim_{N\to\infty}\frac{1}{N}~g_{i,k}(0)=\lambda_{i,k}\alpha_i,
\eqno (3.45)
$$
and (3.39), (3.40) we have
$$
\mathbb{P}-\lim_{N\to\infty}\frac{1}{N}\sum_{i\in\mathscr{
I}_k}g_{i,k}(t)=\varrho_k[1-q(t)]. \eqno (3.46)
$$
Then from (3.45) and (3.46) we obtain
$$
\mathbb{
P}-\lim_{N\to\infty}\frac{1}{N}~g_{i,k}(t)=\lambda_{i,k}\alpha_i[1-q(t)],
\eqno (3.47)
$$
together with
$$
\Phi_s^i(x_k^*)=\frac{\lambda_{i,k}\alpha_k}{\varrho_k}~\Phi_s(x_k^*)
$$
$$~~~~~~~
=\frac{\lambda_{i,k}\alpha_k}{\varrho_k}~q(t).\eqno (3.48)
$$
 Thus, for any $t\ge 0$
$$
\mathbb{P}-\lim_{N\to\infty}\frac{g_{i,j}(t)}{\sum_{l\in\mathscr{
I}_j}g_{l,j}(t)}=\frac{\lambda_{i,k}\alpha_i}{\varrho_k}, \eqno
(3.49)
$$
and for the term of (3.43) under $j=k$ we have
$$
\mathbb{
P}-\lim_{N\to\infty}\frac{1}{N}\int_0^t\frac{g_{i,k}(s)}{\sum_{l\in\mathscr{
I}_j}g_{l,k}(s)}~\mbox{d}\widehat
S_{j}(s)=\frac{\lambda_{i,k}\alpha_i}{\varrho_k}\mu_kt. \eqno
(3.50)
$$
In the following we will use the notation
$$
\beta_{i,j}=\frac{\lambda_{i,j}\alpha_i}{\varrho_j},~~~j=1,2,...,k
\eqno (3.51)
$$
(see also Section 1.4), and then limiting relation (3.50) is
rewritten as
$$
\mathbb{
P}-\lim_{N\to\infty}\frac{1}{N}\int_0^t\frac{g_{i,k}(s)}{\sum_{l\in\mathscr{
I}_j}g_{l,k}(s)}~\mbox{d}\widehat S_{j}(s)=\beta_{i,k}\mu_kt.
\eqno (3.52)
$$

Now, let us consider the case $j\neq k$. Then we have the initial
condition
$$
\mathbb{
P}-\lim_{N\to\infty}\frac{1}{N}~g_{i,j}(0)=\lambda_{i,j}\alpha_i,
\eqno (3.53)
$$
and
$$
\mathbb{P}-\lim_{N\to\infty}\frac{1}{N}\sum_{i\in\mathscr{
I}_j}g_{i,j}(t)=\sum_{i\in\mathscr{I}_j}\lambda_{i,j}\alpha_i
$$
$$~~~~~~~~~~~~~~~~~~~
-[1-q(t)]\sum_{i\in\mathscr{I}_k}\lambda_{i,j}\alpha_i.\eqno
(3.54)
$$
From (3.53) and (3.54) we obtain
$$
\mathbb{
P}-\lim_{N\to\infty}\frac{1}{N}~g_{i,j}(t)=\lambda_{i,j}\alpha_i\{1-I_{i,j}[1-q(t)]\},
\eqno (3.55)
$$
where $I_{i,j}=1$ if $i\in\mathscr{I}_j\bigcap\mathscr{I}_k$, and
$I_{i,j}=0$ otherwise. Thus, for any $t\ge 0$
$$
\mathbb{P}-\lim_{N\to\infty}\frac{g_{i,j}(t)}{\sum_{l\in\mathscr{
I}_j}g_{l,j}(t)}=\frac{\lambda_{i,j}\alpha_i\{1-I_{i,j}[1-q(t)]\}}{\sum_{l\in\mathscr{
I}_j}\lambda_{l,j}\alpha_l\{1-I_{l,j}[1-q(t)]\}}, \eqno (3.56)
$$
and for the term of (3.43) for $j=1,2,...,k-1$ we obtain
$$
\mathbb{
P}-\lim_{N\to\infty}\frac{1}{N}\int_0^t\frac{g_{i,j}(s)}{\sum_{l\in\mathscr{
I}_j}g_{l,j}(s)}~\mbox{d}\widehat S_{j}(s)
$$
$$
= \mathbb{
P}-\lim_{N\to\infty}\frac{1}{N}\int_0^t\frac{\lambda_{i,j}\alpha_i\{1-I_{i,j}[1-q(s)]\}}{\sum_{l\in\mathscr{
I}_j}\lambda_{l,j}\alpha_l\{1-I_{l,j}[1-q(s)]\}}~\mbox{d}\widehat
S_j(s)
$$
$$
= \mathbb{P}-\lim_{N\to\infty}\frac{1}{N}\Big\{\widehat
S_j(t)~\frac{\lambda_{i,j}\alpha_i\{1-I_{i,j}[1-q(t)]\}}{\sum_{l\in\mathscr{
I}_j}\lambda_{l,j}\alpha_l\{1-I_{l,j}[1-q(t)]\}}~~~~~
$$
$$
~~~~~~~~~-\int_0^t\widehat
S_j(s)~\mbox{d}\Big[\frac{\lambda_{i,j}\alpha_i\{1-I_{i,j}[1-q(s)]\}}{\sum_{l\in\mathscr{
I}_j}\lambda_{l,j}\alpha_l\{1-I_{l,j}[1-q(s)]\}}\Big]\Big\}
$$
$$
=
\mu_j\Big\{t~\frac{\lambda_{i,j}\alpha_i\{1-I_{i,j}[1-q(t)]\}}{\sum_{l\in\mathscr{
I}_j}\lambda_{l,j}\alpha_l\{1-I_{l,j}[1-q(t)]\}}~~~~~~~~~~~~~~~~~~~~~~~~
$$
$$
~~~~-\int_0^t
s~\mbox{d}\Big[\frac{\lambda_{i,j}\alpha_i\{1-I_{i,j}[1-q(s)]\}}{\sum_{l\in\mathscr{
I}_j}\lambda_{l,j}\alpha_l\{1-I_{l,j}[1-q(s)]\}}\Big]\Big\}
$$
$$
=\mu_j\int_0^t\frac{\lambda_{i,j}\alpha_i\{1-I_{i,j}[1-q(s)]\}}{\sum_{l\in\mathscr{
I}_j}\lambda_{l,j}\alpha_l\{1-I_{l,j}[1-q(s)]\}}~\mbox{d}s.~~~~~~~~~~~~~~~~~\eqno
(3.57)
$$

In view of (3.52), (3.57) and (3.48), (3.51) system of equations
(3.41) is rewritten as
$$
x_{i,k}^*(t)= \int_0^t\lambda_{i,k}\Big[\alpha_i-
\beta_{i,k}q(t)\Big]\mbox{d}s-\beta_{i,k}\mu_kt,~~~~~~~~~~~~~~
\eqno (3.58)
$$
$$
x_{i,j}^*(t)= \int_0^t\lambda_{i,j}\Big[\alpha_i-
\beta_{i,k}q(t)\Big]\mbox{d}s~~~~~~~~~~~~~~~~~~~~~~~~~~~~~
$$$$~~~~~~~~~~~~~~~~~~~~~~~~~~
-\mu_j\int_0^t\frac{\lambda_{i,j}\alpha_i\{1-I_{i,j}[1-q(s)]\}}{\sum_{l\in\mathscr{
I}_j}\lambda_{l,j}\alpha_l\{1-I_{l,j}[1-q(s)]\}}~\mbox{d}s, \eqno
(3.59)
$$
$$
j=1,2,...,k-1.
$$
From (3.42), (3.58) and (3.59) we obtain the following solution:
$$
x_{i,k}^*(t)=\beta_{i,k}q(t), \eqno (3.60)
$$
$$
x_{i,j}^*(t)=\lambda_{i,j}\alpha_it-\mu_j\int_0^t\frac{\lambda_{i,j}\alpha_i\{1-I_{i,j}[1-q(s)]\}}{\sum_{l\in\mathscr{
I}_j}\lambda_{l,j}\alpha_l\{1-I_{l,j}[1-q(s)]\}}~\mbox{d}s~
$$
$$
-\lambda_{i,j}\alpha_i\beta_{i,k}\int_0^tq(s)\mbox{d}s,\eqno
(3.61)
$$
$$
j=1,2,...,k-1.
$$
From (3.60) and (3.61) we obtain the final solution of the system:
$$
x_k^*(t)=q(t),\eqno (3.62)
$$
and
$$
x_j^*(t)=(\varrho_j\mu_j-\mu_j)t~~~~~~~~~~~~~~~~~~~~~~~~~~~~~~
$$
$$
-\sum_{i\in\mathscr{I}_j\bigcap\mathscr{I}_k}
\lambda_{i,j}\alpha_i\beta_{i,k}\int_0^tq(s)\mbox{d}s, \eqno
(3.63)
$$
$$
j=1,2,...,k-1.
$$

Our next step is to prove the following

\medskip
\noindent {\bf Lemma 3.2.} For any fixed $t>0$ and $\varepsilon>0$
$$
\lim_{N\to\infty}\mathbb{P}\Big\{\sup_{s\le t}\big| x_j(s)-
x_j^*(s)\big|>\varepsilon\Big\}=0,
$$
$$
j=1,2,...,k.
$$

\noindent {\bf Proof.} To prove this lemma only we have to show
that the quadratic characteristics of the square integrable local
martingales $m_{j}(t)$, $j=1,2,...,k$, vanishes in probability,
i.e. for every $\varepsilon>0$
$$
\lim_{N\to\infty}\mathbb{P}\{\big<m_j\big>_t\ge\varepsilon\}=0,~~~\eqno
(3.64)
$$
$$
j=1,2,...,k
$$
(see Kogan and Liptser [16], Lemma 6.1 and Abramov [1], Lemma 2).
Since
$$
\big<M_{j}\big>_t\le\widehat A_j(t)+\widehat S_j(t)
$$
$$~~~~~~~~~~~~~~~~
\le t\sum_{i\in\mathscr{I}_j}\lambda_{i,j}N_i+\widehat
S_j(t),\eqno (3.65)
$$
then taking into account (3.15) and
$$
\big<m_{j}\big>_t\le t\sum_{i\in\mathscr{I}_j}
\frac{\lambda_{i,j}}{N_i}+\frac{1}{N^2}\widehat S_j(t) \eqno
(3.66)
$$
we obtain the desired statement of Lemma 3.2. The lemma is
proved.\\

Applying Liptschitz conditions (3.26), (3.27), for any
$\varepsilon>0$ we have
$$
\lim_{N\to\infty}\mathbb{P}\big\{\sup_{s\le
t}q_j(s)\ge\varepsilon\big\}=0, \eqno (3.67)
$$
$$
j=1,2,...,k-1,
$$
and
$$
\lim_{N\to\infty}\mathbb{P}\big\{\sup_{s\le
t}\big|q_k(s)-x_k^*(s)\big|\ge\varepsilon\big\}=0. \eqno (3.68)
$$

\section{A stability of the queue-length processes in the client
stations} \noindent Now, let us study a question on stability,
i.e. existence of the limiting generalized functions of the
non-stationary probabilities
$$
\lim_{N\to\infty}\int_0^t\mathbb{P}\{Q_j(s)=l\}\mbox{d}s, \
l=0,1,... \eqno(4.1)
$$
for all $t>0$ at the client stations $j=1,2,...,k-1$, satisfying
$$
\lim_{N\to\infty}\sum_{l=0}^\infty\int_0^t\mathbb{P}\{Q_j(s)=l\}\mbox{d}s=t.
$$

The existence of the limiting generalized functions does not mean
existence of the limiting non-stationary probabilities
$$
\lim_{N\to\infty}\mathbb{P}\{Q_j(s)=l\}, \ l=0,1,... \eqno(4.2)
$$

The proof of the lemma below is based on the same idea that the
corresponding proof in Abramov [1] and partially repeats that
proof. However, there is a place in the proof of stability in
Abramov [1] that should be reconsidered and improved, and the
proof given below is doing that. Furthermore, the fact, that there
are several server stations with different behavior of the
queue-lengths, adds a number of significant features as well.

\medskip
\noindent {\bf Lemma 4.1.} Under the assumptions given in the
paper for all $j=1,2,...,k-1$ there exist the stationary
queue-length processes $Q_j^*(s)$, $j=1,2,...,k-1$, satisfying the
inequalities
$$
\liminf_{N\to\infty}\mathbb{P}\{Q_j(t)=
l\}\le\mathbb{P}\{Q_j^*(s)=
l\}\le\limsup_{N\to\infty}\mathbb{P}\{Q_j(t)=l\}
$$
$$
l=0,1,...
$$


\medskip
\noindent {\bf Proof.} To prove this lemma let us denote the
queue-length process in the server station $i$ by $\Sigma_i(t)$.
As earlier, $\mathscr{J}_i$ denote a set of indexes $j$ where
$\lambda_{i,j}>0$. We have
$$
\Sigma_i(t)=N_i-\sum_{l\in\mathscr{J}_l}Q_{i,l}(t),~~~i=1,2,...,r.
\eqno (4.3)
$$
Normalization of (4.3) yields
$$
\lim_{N\to\infty}\frac{1}{N}\Sigma_i(t)=\alpha_i-\sum_{l\in\mathscr{
J}_i}q_{i,l}(t),~~~i=1,2,...,r. \eqno (4.4)
$$
Then, according to (3.67), (3.68) and (3.62), (3.63) we have the
following two types of nodes. If $i$ does not belong to
$\mathscr{I}_k$ then
$$
\mathbb{P} - \lim_{N\to\infty}\frac{1}{N}\Sigma_i(t)=\alpha_i,
\eqno (4.5)
$$
for any $t>0$, otherwise, if $i\in\mathscr{I}_k$ then for all
$t>0$
$$
\mathbb{P} - \lim_{N\to\infty}\frac{1}{N}\Sigma_i(t)<\alpha_i.
\eqno (4.6)
$$
More exactly, in the last case
$$
\mathbb{P} -
\lim_{N\to\infty}\frac{1}{N}\Sigma_i(t)=\alpha_i[1-q(t)\beta_{i,k}].
\eqno (4.7)
$$
The relation (4.4) means that
$$
\lim_{N\to\infty}\mathbb{P}\Big\{\sup_{s\le
t}\frac{1}{N}\Sigma_i(s)\le\alpha_i\Big\}=1,\eqno (4.8)
$$
$$
i=1,2,...,r.
$$

Prove that
$$
\mathbb{P}\Big\{\lim_{N\to\infty}\sup_{s\le
t}[A_j(s)-S_j(s)]<\infty\Big\}=1,\eqno (4.9)
$$
$$
j=1,2,...,k-1.
$$
Indeed, it follows from (4.8) that $A_{j}(t)\le\pi_j(t)$ where
$\pi_j(t)$ is a Poisson process with parameter
$\sum_{i\in\mathscr{I}_j}\lambda_{i,j}N_i$ given on the same
probability space as the process $A_{j}(t)$. Therefore the problem
is reduced to prove
$$
\mathbb{P}\Big\{\lim_{N\to\infty}\sup_{s\le
t}[\pi_j(s)-S_j(s)]<\infty\Big\}=1,
 \eqno (4.10)
$$
$$
j=1,2,...,k-1.
$$
Obviously, that $\mathbb{P}$-a.s.
$$
\lim_{N\to\infty}\frac{\pi_j(t)-S_j(t)}{Nt}=\varrho_j\mu_j-\mu_j<1.
\eqno (4.11)
$$
Therefore,
$$
\lim_{N\to\infty}\big[\pi_j(t)-S_j(t)\big]=-\infty, \eqno (4.12)
$$
and (4.10), (4.9) follow by the fact that both $\pi_j(t)$ and
$S_j(t)$ are c\'adl\'ag functions. Next, let $ \widetilde
Q_j(t)=\pi_j(t)-D_j(t) $ where $D_j(t)$ is a function defined in
(2.1) and (2.3). According to the Skorokhod reflection principle
$$
\widetilde Q_j(t)=\pi_j(t)-S_j(t)-\inf_{s\le
t}\big[\pi_j(t)-S_j(t)\big]. \eqno (4.13)
$$
Because of the strict stationarity and ergodicity of increments of
the process $\pi_j(t)-S_j(t)$, from (4.13) we obtain
$$
\widetilde Q_j(t)=\sup_{s\le t}\big[\big( \pi_j(t)-S_j(t)\big)-
\big(\pi_j(s)-S_j(s)\big)\big], \eqno (4.14)
$$
and in addition
$$
\sup_{s\le t}\big[\big( \pi_j(t)-S_j(t)\big)-
\big(\pi_j(s)-S_j(s)\big)\big]{\buildrel \mbox{d}\over
=}\sup_{s\le t}\big[\pi_j(s)-S_j(s)\big]. \eqno (4.15)
$$
Thus we have
$$
\widetilde Q_j(t) {\buildrel \mbox{d}\over =}\sup_{s\le
t}\big[\pi_j(s)-S_j(s)\big]. \eqno (4.16)
$$
Taking into account that $A_j(t)\le \pi_j(t)$ we obtain
$$
\liminf_{N\to\infty}\mathbb{P}\{Q_j(t)=l\}\ge\mathbb{P}\{
\widetilde Q_j(t)=l\}
$$
$$
~~~~~~~~~~~~~~~~~~~~~~~~~~~~~~~~~~~~~~~~~~~=\lim_{N\to\infty}\mathbb{
P}\big\{\sup_{s\le t}\big[\pi_j(s)-S_j(s)\big]=l\big\}.\eqno
(4.17)
$$

Next, let us introduce a Poisson process $\Pi_j(z)$ with parameter
$$
\sum_{i\in\mathscr{I}_j}\lambda_{i,j}N_i[1-\beta_{i,k}q(t)] \eqno
(4.18)
$$
Assuming that both processes $A_j(z)$ and $\Pi_j(z)$ are given on
the same probability space we have $A_j(z)>\Pi_j(z)$ for all
$z>0$, and because of
$$
\sum_{i\in\mathscr{
I}_j}\lambda_{i,j}N_i[1-\beta_{i,k}q(t)]<\mu_j, \eqno (4.19)
$$
then analogously to (4.17) we have the following
$$
\limsup_{N\to\infty}\mathbb{P}\{Q_j(t)=l\}\le
\lim_{N\to\infty}\mathbb{ P}\big\{\sup_{s\le
t}\big[\Pi_j(s)-S_j(s)\big]=l\big\}. \eqno (4.20)
$$
Inequalities (4.17) and (4.20) together with the results of
convergence in probability (see ref. (3.67), (3.68) and Lemma 3.2)
allow us to conclude that there exists a stationary process
$Q_j^*(t)$, and for some $\alpha(t)$, $0\le\alpha(t)\le 1$, we
have
$$
\mathbb{
P}\{Q_j^*(t)=l\}=\alpha(t)\liminf_{N\to\infty}\mathbb{P}\{Q_j(t)=l\}
$$
$$~~~~~~~~~~~~~~~~~~~~~~~~~~~~~~~~
+[1-\alpha(t)]\limsup_{N\to\infty}\mathbb{P}\{Q_j(t)=l\},\eqno
(4.21)
$$
and the lemma is therefore proved.\\

Notice, that existence of some stationary process $Q_j^*(t)$ does
not mean that there exist limiting non-stationary queue-length
distributions (4.2). However, then there exists the limiting
generalized queue-length distributions (4.1). Notice also that in
the case when $S_j(t)$, $j=1,2,...,k$, are Poisson processes, the
limiting non-stationary queue-length distributions (4.2) do exist
and coincide with the distribution of the processes $Q_j^*(t)$.
The existence of limiting non-stationary distributions in the last
case follows from the Chapman-Kolmogorov equations, which can be
written in explicit form.

\section{The proof of the main result and special cases}
\noindent We start from the proof of the main result.

\medskip \noindent {\bf Proof of Theorem 1.1.} The proof of this
theorem is analogous to the corresponding proof of Theorem 3 of
Abramov [1]. Therefore we pay more attention to the new features
where detailed explanation is necessary.

First, we have the representation:
$$
\lim_{N\to\infty}\frac{1}{\mu_jN}\mathbb{E}\int_0^t{\bf
I}\{Q_j(s-)=l\}\mbox{d}S_j(s)
$$
$$
=\lim_{N\to\infty}\int_0^t\mathbb{P}\{Q_j[S_j^*(s)]=l\}\mbox{d}s,
\eqno (5.1)
$$
$$
~~~l=0,1,...,
$$
where $S_j^*(s)$ is introduced in Section 1.4. (For the proof of
(5.1) see Appendix A.)

Then, introducing
$$
A_j^*(t)=\inf\{s>0: A_j(s)=A_j(t)\},\eqno (5.2)
$$
$$
j=1,2,...,k-1
$$
we have the representation analogous to (5.1)
$$
\lim_{N\to\infty}\frac{1}{\sum_{i\in\mathscr{I}_j}\lambda_{i,j}N_i}
\mathbb{E}\int_0^t{\bf I}\{Q_j(s-)=l\}\mbox{d}A_j(s)
$$
$$
=\lim_{N\to\infty}\int_0^t\mathbb{P}\{Q_j[A_j^*(s)]=l\}\mbox{d}s,\eqno
(5.3)
$$
$$
~~~l=0,1,...~.
$$
(The proof of (5.3) is completely analogous to the proof of (5.1)
given in Appendix A.)

Next, for all $l=1,2,...$;~ $j=1,2,...,k-1$ and $t>0$ we have the
relation connecting the number of up- and down-crossings:
$$
\sum_{u=1}^{A_j(t)}{\bf
I}\{Q_j(\tau_{j,u}-)=l-1\}=\sum_{u=1}^{S_j(t)}{\bf
I}\{Q_j(\sigma_{j,u}-)=l\}
$$
$$~~~~~~~~~~~~~~~~~~~~~~~
+{\bf I}\{Q_j(t)\ge l\}, \eqno (5.4)
$$
where $\tau_{j,1}$, $\tau_{j,2},...$ are the moments of arrival to
node $j$, and $\sigma_{j,1}$, $\sigma_{j,2},...$ are the moments
of departure from node $j$.

It follows from (5.1), (5.3) and (5.4) that
$$
\lim_{N\to\infty}\frac{1}{N}\int_0^t{\bf
I}\{Q_j(s-)=l-1\}\mbox{d}A_j(s)
$$
$$
= \lim_{N\to\infty}\frac{1}{N}\int_0^t{\bf
I}\{Q_j(s-)=l\}\mbox{d}S_j(s),\eqno (5.5)
$$
$$
l=1,2,...~.
$$
Taking into account (3.20) the left-hand side of (5.5) can be
rewritten as
$$
\lim_{N\to\infty}\frac{1}{N}\mathbb{E}\int_0^t{\bf
I}\{Q_j(s-)=l-1\}\mbox{d}A_j(s)
$$
$$
=\lim_{N\to\infty}\frac{1}{N}\mathbb{E}\int_0^t{\bf
I}\{Q_j(s-)=l-1\}\mbox{d}\widehat A_j(s). \eqno (5.6)
$$
(The technical details of the proof of that see in Appendix B.)

From (2.13), (2.15) and asymptotic results of Section 3 (see
(3.39), (3.40) and (3.48) together with notation (3.51)) we obtain
$$
\mathbb{P}-\lim_{N\to\infty}\frac{1}{N}\widehat
A_j(t)=\varrho_j\mu_j \Big[1-q(t)\sum_{i\in\mathscr{I}_j\bigcap
\mathscr{I}_k}\beta_{i,j}\beta_{i,k}\Big],\eqno (5.7)
$$
$$
~~~ j=1,2,...,k-1.
$$
Therefore, substituting (5.7) for (5.6) we obtain
$$
\lim_{N\to\infty}\frac{1}{N}\int_0^t{\bf
I}\{Q_j(s-)=l-1\}\mbox{d}\widehat A_j(s)
$$
$$=
\varrho_j\mu_j\int_0^t\mathbb{P}\{Q_j(s)=l-1\}
\Big[1-q(s)\sum_{i\in\mathscr{I}_j\bigcap \mathscr{
I}_k}\beta_{i,j}\beta_{i,k}\Big]\mbox{d}s.\eqno (5.8)
$$
In turn, substituting (5.8) for (5.6) and (5.5) we obtain
$$
\lim_{N\to\infty}\frac{1}{N}\int_0^t{\bf
I}\{Q_j(s-)=l\}\mbox{d}S_j(s)$$$$= \varrho_j\mu_j\int_0^t\mathbb{
P}\{Q_j(s)=l-1\}\Big[1-q(s)\sum_{i\in\mathscr{ I}_j\bigcap
\mathscr{I}_k}\beta_{i,j}\beta_{i,k}\Big]\mbox{d}s.\eqno (5.9)
$$
Then, keeping in mind (5.1) proves the theorem. The theorem is
proved.\\

Let us now prove the special cases of this theorem.\\

\noindent{\bf Proof of Corollary 1.2}. Notice that
$$
\lim_{N\to\infty}\frac{1}{N}\mathbb{E}\int_0^t{\bf
I}\{Q_j(s-)=l\}\mbox{d}S_j(s)
$$
$$=
\lim_{N\to\infty}\frac{1}{N}\mathbb{E}\int_0^t{\bf
I}\{Q_j(s-)=l\}\mbox{d}\widehat S_j(s),\eqno (5.10)
$$
$$
l=0,1,...~.
$$
(For the proof see Appendix B by replacing there the processes
$A_j(t)$ and $\widehat A_j(t)$ with the corresponding processes
$S_j(t)$ and $\widehat S_j(t)$.) Then limiting relation (5.1) can
be rewritten as
$$
\lim_{N\to\infty}\frac{1}{\mu_jN}\mathbb{E}\int_0^t{\bf
I}\{Q_j(s-)=l\}\mbox{d}\widehat S_j(s)
$$
$$
=\lim_{N\to\infty}\int_0^t\mathbb{P}\{Q_j[S_j^*(s)]=l\}\mbox{d}s,
\eqno (5.11)
$$
$$
~~~l=0,1,...~.
$$
According to assumption $S_j(t)$ is the Poisson process with rate
$\mu_jN$, and therefore, $\widehat S_j(t)=\mu_jNt$. Then, from
(5.11) we obtain
$$
\lim_{N\to\infty}\mathbb{P}\{Q_j[S_j^*(t)]=l\}=\lim_{N\to\infty}\mathbb{
P}\{Q_j(t)=l\},\eqno (5.12)
$$
$$
l=0,1,...~,
$$
and the result follows.\\

\noindent {\bf Proof of Corollary 1.3.} Indeed, in this case
$\beta_{1,j}=\beta_{1,k}=1$ for
all $j=1,2,...,k-1$, and the result follows from Theorem 1.1.\\

\noindent {\bf Proof of Corollary 1.4.} Indeed, in this case both
$\beta_{i,k}$ and $\beta_{i,j}$ are positive for the same set of
indexes $i$ and therefore
$$
\sum_{i\in\mathscr{I}_j\bigcap\mathscr{I}_k}\beta_{i,k}\beta_{i,j}=1,
\eqno (5.13)
$$
and the result follows.\\

\noindent
{\bf Proof of Corollary 1.5.} The proof follows immediately from Theorem 1.1.\\

\noindent {\bf Proof of Corollary 1.6.} Indeed, when $\varrho_k=1$
we obtain $q(t)=0$, and the result follows from Theorem 1.1.

\section*{APPENDIX A}

\noindent {\bf Proof of relation (5.1).} Take a small interval
$\mathscr{U}=(u, u+\mbox{d}u]$. Denote
$$
n_1=\min\{ n: \sigma_{j,n}\in \mathscr{U}\}
$$
$$
n_2=\max\{ n: \sigma_{j,n}\in \mathscr{U}\},
$$
where the notation for $\sigma_{j,n}$ is given in Section 1.4.
Recall that $\sigma_{j,n}=\sum_{i=1}^n\xi_{j,i}$, where
$\{\xi_{j,i}\}_{i\ge 1}$ are the increments associated with the
point process $S_j(t)$. Then, according to Lemma 3.1 we obtain
that $N^{-1}[S_j(u+\mbox{d}u)-S_j(u)]{\buildrel\mathbb{
P}\over\to}\mu_j\mbox{d}u$.

Let us now apply the Cesaro theorem (e.g. P\'olya and Szeg\"o
[30]): {\it if a sequence $\{a_n\}$ converges to $S$, then also a
sequence $S_N=N^{-1}\sum_{i=1}^Na_i$ converges to $S$.} Applying
this result and the Lebesgue theorem on dominated convergence we
obtain:
$$
\lim_{N\to\infty}\mathbb{P}\{Q_j(S^*(u+\mbox{d}u)=l)
$$
$$
=\mathbb{
P}-\lim_{N\to\infty}\frac{1}{S_j(u+\mbox{d}u)-S_j(u)}\sum_{i=n_1}^{n_2}\mathbb{
P}\{Q_j(\sigma_{j,i}-)=l\}
$$
$$
=\lim_{N\to\infty}\frac{1}{\mu_jN\mbox{d}u}\mathbb{
E}\sum_{i=n_1}^{n_2}\mathbb{P}\{Q_j(\sigma_{j,i}-)=l\}
$$
$$
=\lim_{N\to\infty}\frac{1}{\mu_jN\mbox{d}u}\sum_{\alpha_1\le\alpha_2}\mathbb{
E}\Big\{\sum_{i=\alpha_1}^{\alpha_2} \mathbb{
P}\{Q_j(\sigma_{j,i}-)=l\}\Big|n_1=\alpha_1, n_2=\alpha_2\Big\}
$$
$$\times
\mathbb{P}\{n_1=\alpha_1, n_2=\alpha_2\}
$$
$$
=\lim_{N\to\infty}\frac{1}{\mu_jN\mbox{d}u}\sum_{\alpha_1\le\alpha_2}\mathbb{
E}\Big\{\sum_{i=\alpha_1}^{\alpha_2}{\bf
I}\{Q_j(\sigma_{j,i}-)=l\}\Big|n_1=\alpha_1, n_2=\alpha_2\Big\}
$$
$$\times
\mathbb{P}\{n_1=\alpha_1, n_2=\alpha_2\}
$$
$$
=\lim_{N\to\infty}\frac{1}{\mu_jN\mbox{d}u}\mathbb{
E}\sum_{i=n_1}^{n_2}{\bf I}\{Q_j(\sigma_{j,i}-)=l\}
$$
$$
=\frac{1}{\mbox{d}u}\lim_{N\to\infty}\frac{1}{\mu_jN}\mathbb{E}({\bf
I}\{Q_j(u^*-)=l\}[S_j(u+\mbox{d}u)-S_j(u)]),
$$
where $u^*\in\mathscr{U}$. Relation (5.1) follows.

\section*{APPENDIX B}

\noindent{\bf Proof of relation (5.6)}. Rewrite the left-hand side
of (5.6) as
$$~~~~~
\lim_{N\to\infty}\frac{1}{N}\mathbb{E}\int_0^t{\bf
I}\{Q_j(s-)=l-1\}\mbox{d}A_j(s)
$$
$$
=\lim_{N\to\infty}\frac{1}{N}\mathbb{E}\int_0^t{\bf
I}\{Q_j(s-)=l-1\}\mbox{d}\widehat A_j(s)
$$
$$~~~~~~~~~~~
+ \lim_{N\to\infty}\frac{1}{N}\mathbb{E}\int_0^t{\bf
I}\{Q_j(s-)=l-1\}\mbox{d}[A_j(s)-\widehat A_j(s)]\eqno (B.1)
$$
It follows from (3.20) that, as $N\to\infty$, the term
$$
\frac{1}{N}[A_j(t)-\widehat A_j(t)] \eqno (B.2)
$$
vanishes in probability, and also
$$
\Big|\mathbb{E}\int_0^t{\bf
I}\{Q_j(s-)=l-1\}\mbox{d}[A_j(s)-\widehat A_j(s)]\Big|
$$
$$
 \le |\mathbb{E}(A_j(t)-\widehat A_j(t))|. \eqno (B.3)
$$
Therefore, as $N\to\infty$, from (B.2), (B.3) and the Lebesgue
theorem on dominated convergence we obtain
$$
\lim_{N\to\infty}\frac{1}{N}\mathbb{E}\int_0^t{\bf
I}\{Q_j(s-)=l-1\}\mbox{d}[A_j(s)-\widehat A_j(s)]=0,
$$
and (5.6) follows.

\section*{Acknowledgement}

The author thanks the anonymous referees for their comments
helping substantially improve the paper.

\section*{References}

\indent

[1] V.M.Abramov, A large closed queueing network with autonomous
service and bottleneck. Queueing Systems Theory Appl. 35 (2000)
23-54.

\smallskip

[2] V.M.Abramov, Some results for large closed queueing networks
with and without bottleneck: Up- and down-crossings approach.
Queueing Systems Theory Appl. 38 (2001) 149-184.

\smallskip

[3] S.V.Anulova and R.Sh.Liptser, Diffusion approximation for
processes with normal reflection. Theory Prob. Appl. 35 (1990)
413-423.

\smallskip

[4] A.A.Borovkov, {\em Stochastic Processes in Queueing Theory.}
(Springer-Verlag, Berlin, 1976).

\smallskip

[5] A.A.Borovkov, {\em Asymptotic Methods in Queueing Theory.}
(John Wiley, New York, 1984).

\smallskip

[6] Yi-Ju Chao, Weak convergence of a sequence of semimartingales
to a diffusion with discontinuous drift and diffusion
coefficients. Queueing Systems Theory Appl. 42 (2002)  153-188.

\smallskip

[7] H.Chen and A.Mandelbaum, Discrete flow networks: Bottleneck
analysis and fluid approximations. Math. Operat. Res. 16 (1991)
408-446.

\smallskip

[8] H.Chen and A.Mandelbaum, Discrete flow networks: Diffusion
approximations and bottlenecks. Ann. Prob. 19 (1991) 1463-1519.

\smallskip

[9] C.Dellacherie, {\em Capacit\'es et Processus Stochastiques}.
(Springer-Verlag, Berlin, 1972).

\smallskip

[10] C.Knessl,  B.J.Matkovsky, Z.Schuss and C.Tier, Asymptotic
analysis of a state dependent $M/G/1$ queueing system. SIAM J.
Appl. Math. 46(3) (1986) 483-505.

\smallskip

[11] C.Knessl,  B.J.Matkovsky, Z.Schuss and C.Tier, The two
repairmen problem: A finite source $M/G/2$ queue. SIAM J. Appl.
Math. 47(2) (1987) 367-397.

\smallskip

[12] C.Knessl,  B.J.Matkovsky, Z.Schuss and C.Tier, A Markov
modulated $M/G/1$ queue. I. Stationary distribution. Queueing
Systems Theory Appl. 1 (1987) 355-374.

\smallskip

[13] C.Knessl,  B.J.Matkovsky, Z.Schuss and C.Tier, A Markov
modulated $M/G/1$ queue. II. Busy period and time for buffer
overflow. Queueing Systems Theory Appl. 1 (1987) 375-399.

\smallskip

[14] C.Knessl,  B.J.Matkovsky, Z.Schuss and C.Tier,
 Busy period distribution in state dependent queues.
Queueing Systems Theory Appl. 2 (1987) 285-305.

\smallskip

[15] C.Knessl and C.Tier, Asymptotic expansions for large closed
queueing networks. J. Assoc. Comput. Mach. 37 (1990) 144-174.

\smallskip

[16] Y.Kogan, and R.Sh.Liptser, Limit non-stationary behavior of
large closed queueing networks with bottlenecks. Queueing Systems
Theory Appl. 14 (1993) 33-55.

\smallskip

[17] Y.Kogan, R.Sh.Liptser and M.Shenfild, State dependent Bene\^s
buffer model with fast loading and output rates. Ann. Appl.
Probab. 5 (1995) 97-120.

\smallskip

[18] Y.Kogan, R.Sh.Liptser and A.V.Smorodinskii, Gaussian
diffusion approximation of closed Markov model of computer
networks. Probl. Inform. Transmission 22 (1986) 38-51.

\smallskip

[19] T.Konstantopoulos, S.N.Papadakis and J.Walrand, Functional
approximation theorems for controlled queueing systems. J. Appl.
Prob. 31 (1994) 765-776.

\smallskip

[20] E.V.Krichagina, R.Sh.Liptser and A.A.Puhalskii, Diffusion
approximation for system with arrivals depending on queue and
arbitrary service distribution. Theory Probab. Appl. 33 (1988)
114-124.

\smallskip

[21] E.V.Krichagina and A.A.Puhalskii, A heavy traffic analysis of
closed queueing system with $GI/\infty$ service senter. Queueing
Systems Theory Appl. 25 (1997) 235-280.

\smallskip

[22] N.V.Krylov and R.Liptser, On diffusion approximation with
discontinuous coefficients. Stoch. Proc. Appl. 102 (2002) 235-264.

\smallskip

[23] R.Sh.Liptser, A large deviation problem for simple queueing
model. Queueing Systems Theory Appl. 14 (1993) 1-32.

\smallskip

[24] R.Sh.Liptser and A.N.Shiryayev, {\em Statistics of Random
Processes}, Vols. I, II. (Springer-Verlag, Berlin, 1977-1978).

\smallskip

[25] R.Sh.Liptser and A.N.Shiryayev, {\em Theory of Martingales}.
(Kluwer, Dordrecht, 1989).

\smallskip

[26] A.Mandelbaum and W.A.Massey, Strong approximations for time
dependent queues. Math. Operat. Res. 20 (1995) 33-64.

\smallskip

[27] A.Mandelbaum, W.A.Massey and M.I.Reiman, Strong
approximations for Markovian service networks. Queueing Systems
Theory Appl. 30 (1998) 149-201.

\smallskip
[28] A.Mandelbaum and G.Pats, State-dependent queues:
approximations and applications. In: {\em IMA Volumes in
Mathematics and Its Applications.} (F.P.Kelly and R.J.Williams
eds), 71, pp. 239-282. (Springer-Verlag, Berlin, 1995).

\smallskip

[29] A.Mandelbaum and G.Pats, State-dependent stochastic networks.
Part I. Approximations and applications with continuous diffusion
limits. Ann. Appl. Probab. 8 (1998) 569-646.

\smallskip
[30] G.P\'olya and G.Szeg\"o, {\em Aufgaben und Lehrsatze aus der
Analysis. Erster Band: Reihen, Integralrechnung,
Functionentheorie.} (Springer-Verlag, Berlin, 1964).

\smallskip

[31] M.I.Reiman and B.Simon, A network of priority queues in heavy
traffic: One bottleneck station. Queueing Systems Theory Appl. 6
(1990) 33-58.

\smallskip

[32] A.V.Skorokhod, Stochastic equations for diffusion processes
in a bounded region. Theory Probab. Appl. 6 (1961) 264-274.

\smallskip

[33] H.Tanaka, Stochastic differential equations with reflected
boundary condition in convex regions. Hiroshima Math. J. 9 (1979)
163-177.

\smallskip

[34] W.Whitt, Open and closed models for networks of queues. AT\&T
Bell. Lab. Tech. J. 63 (1984) 1911-1979.

\smallskip

[35] R.J.Williams, On approximation of queueing networks in heavy
traffic. In: {\em Stochastic Networks. Theory and Application},
eds. F.P.Kelly, S.Zachary, I.Ziedins (Oxford Univ. Press, Oxford,
1996), pp. 35-56.

\smallskip

[36] R.J.Williams, An invariance principle for semimartingale
reflecting Brownian motion in orthant. Queueing Systems Theory
Appl. 30 (1998) 5-25.

\smallskip
[37] R.J.Williams, Diffusion approximation for open multiclass
queueing networks: Sufficient conditions involving state space
collapse. Queueing Systems Theory Appl. 30 (1998) 27-88.

\end{article}

\end{document}